\documentclass[11pt]{article}
\usepackage{times}
\usepackage{amsmath,amsfonts,amssymb,latexsym,epsfig,a4}
\usepackage{color,epsf}
\usepackage{graphicx}

\usepackage{alltt}
\usepackage{url}

\newcommand{\e}{\ensuremath{\mathrm{e}}}
\newcommand{\ad}{\ensuremath{\mathrm{ad}}}

\begin{document}
\title{Magnus integrators for linear and quasilinear delay differential equations
}

\author{Ana Arnal\footnote{ Email: \texttt{ana.arnal@uji.es}. ORCID: 0000-0002-3283-3379}
 \and
 Fernando Casas\footnote{Email: \texttt{Fernando.Casas@uji.es}. ORCID: 0000-0002-6445-279X} 
 \and Cristina Chiralt\footnote{Email: \texttt{chiralt@uji.es}. ORCID: 0000-0003-0925-4034}
}


%
\maketitle

\begin{center}
Institut de Matem\`atiques i Aplicacions de Castell\'o (IMAC) and  De\-par\-ta\-ment de
Ma\-te\-m\`a\-ti\-ques, Universitat Jaume I,
  E-12071 Cas\-te\-ll\'on, Spain.
\end{center}

\vspace*{0.5cm}

\begin{abstract}
A procedure to numerically integrate non-autonomous linear delay differential equations is presented. It is based on the use of an spectral discretization
of the delayed part to transform the original problem into a matrix linear ordinary differential equation which is subsequently solved with numerical
integrators obtained from the Magnus expansion. The algorithm can be used in the periodic case
to get both accurate approximations of the characteristic multipliers and the 
solution itself. In addition, it can be extended to deal with certain quasilinear delay equations.

\


\noindent

\end{abstract}\bigskip



\section{Introduction}

We are primarily interested in the numerical integration of a system of linear delay differential equations (DDEs) with a single 
 discrete time delay of the form
\begin{equation} \label{dde.1}
\aligned
  & \frac{dx}{dt} = A(t) \, x(t) + B(t) \, x(t - \tau), \quad t \ge 0 \\
  & x(t) = \phi(t), \qquad -\tau \le t \le 0,
\endaligned
\end{equation}  
where $x \in \mathbb{R}^d$, $A$ and $B$ are, in general, time-dependent $d \times d$ matrices and $\phi(t)$ is the initial function. Two cases of special
interest arise: when $A$ and $B$ are constant, and when they are periodic of period $T$. Our purpose is to present an algorithm that allows one to
analyze the stability of the problem at $t=T$ in the periodic case  and also get approximations to the solution of (\ref{dde.1}) for arbitrary $t > 0$ in the
general (not necessarily periodic) time dependent case.
The formalism can be readily generalized to problems with $k$ distinct delays $0 < \tau_1 < \cdots < \tau_k$,
\begin{equation} \label{dde.1b}
\aligned
  & \frac{dx}{dt} = A(t) \, x(t) + \sum_{j=1}^k B_j(t) x(t - \tau_j), \quad t \ge 0 \\
  & x(t) = \phi(t), \qquad -\tau \le t \le 0.
\endaligned
\end{equation}  
Moreover, the procedure  is also extended to quasilinear delay equations of the form
\begin{equation} \label{dde.2}
\aligned
  & \frac{dx}{dt} = A \big(x(t-\tau) \big) \, x(t), \quad t \ge 0 \\
  & x(t) = \phi(t), \qquad -\tau \le t \le 0.
\endaligned
\end{equation}  
Equations (\ref{dde.1}) and (\ref{dde.1b}) appear frequently in applications, either as a model of some physical problem or as a tool to analyze its stability when
the time evolution of the unknown variable depends not only on the actual state but also on its past values (see e.g. \cite{breda15sol,insperger11sdf} and
references therein). On the other hand, equation (\ref{dde.2}) has been used to describe SIR-type epidemic models taking into account the latent period, i.e., 
the time when an individual is infected but is not infective \cite{csomos20mti}.


Given the relevance of problems (\ref{dde.1})-(\ref{dde.2}), many numerical procedures have been designed over the years for obtaining approximate
solutions (see \cite{bellen03nmf} and references therein). Among them, the approximation technique consisting in first converting the DDE (\ref{dde.1}) when
$A$ and $B$ are constant into an abstract Cauchy problem \cite{bellen00nso,maset99asi} 
is particularly appealing: essentially, it allows one to discretize the corresponding operator
and then solve numerically the resulting system of ordinary differential equations (ODEs) by standard methods. This procedure has been used to analyze the
stability of linear DDEs, both autonomous \cite{breda15sol} and explicitly time-dependent \cite{butcher11otc}, 
in combination with Chebyshev spectral collocation methods. It is called ``continuous time approximation" in \cite{butcher11otc,sun09amo} and, in particular, provides
spectral accuracy in the determination of the characteristic roots of the system.

We pursue here the same strategy and combine it with the application of numerical integrators based on the Magnus expansion to carry out the
time integration of the resulting non-autonomous 
system of ODEs. We show that this procedure provides more accurate approximations 
than those reported in \cite{breda15sol} for the determination of the characteristic multipliers of (\ref{dde.1b}) with the same number of discretization points. 
Even for the quasilinear case (\ref{dde.2}) it leads to
higher order approximations to the solution than previous schemes also based on the Magnus expansion \cite{csomos20mti}, whereas still preserving 
its qualitative properties.

The plan of the paper is the following. 
In section \ref{sec.2} we briefly summarize the continuous time approximation technique for dealing with linear
DDEs, whereas the main features of the Magnus expansion and some numerical integrators based on it are reviewed in section \ref{sec.3}. The time
integration algorithm is illustrated on several numerical examples in section \ref{sec.4} and the technique is extended in section \ref{sec.5} to quasilinear
problems, and in particular to an epidemic model with delay.
Finally, section \ref{sec.6} contains some concluding remarks.

\section{Continuous time approximation}
\label{sec.2}

\subsection{Linear DDEs as abstract Cauchy problems}

Let us denote by $X$ the state space of continuous functions $C([-\tau,0], \mathbb{R}^d)$, which is a Banach space with the norm
$\|\phi\|_X \equiv \max_{-\tau \le \theta \le 0} \|\phi(\theta)\|_{\infty}$ \cite{breda15sol}. If $x_t \in X$ is the state at time $t$, defined as \cite{hale93itf}
\begin{equation} \label{dde.2b}
  x_t(\theta) \equiv x(t +\theta), \quad \theta \in [-\tau, 0],
\end{equation}
then the linear equation (\ref{dde.1}) in the autonomous case can be written as
\begin{equation} \label{dde.3}
\aligned
  & \dot{x} = L \,  x_t, \qquad t \in \mathbb{R} \\
  & x_0 = \phi \in X,
\endaligned  
\end{equation}
where it is assumed that the initial time $t_0 = 0$,
the dot denotes the right-hand derivative and $L : X \longrightarrow \mathbb{R}^d$ acts on $x_t$ as
\begin{equation} \label{dde.4}
     L \, x_t = A \, x(t) + B \, x(t - \tau).
\end{equation}
In that case, for every
$(t_0, \phi) \in \mathbb{R} \times X$, there exists a unique solution of (\ref{dde.3}) on $[- \tau, +\infty)$, denoted by
$x(t; \phi)$ \cite{breda15sol}:
\[
  x(t; \phi) = \left\{ \begin{array}{ll}
     \phi(0) + \displaystyle \int_0^t L \, x_s \, ds, & \quad t \ge 0 \\
     \phi(t), & \quad t \in [-\tau, 0].
     \end{array} \right.
\]     
It is then possible to apply the theory of one-parameter strongly continuous semigroups (also called
$C_0$-semigroups) in this setting. More specifically, it has been shown that the family $\{ U(t)\}_{t \ge 0}$ 
of operators $U(t): X \longrightarrow X$
associating to the initial function $\phi$ the state $x_t$ at time $t \ge 0$, i.e, $U(t) \phi = x_t(\cdot; \phi)$, defines a 
$C_0$-semigroup of linear and bounded operators on $X$ whose infinitesimal generator is the linear unbounded operator
$\mathcal{A}: \mathcal{D}(\mathcal{A})  \subseteq X \longrightarrow X$ given by \cite{engel00ops}
\begin{equation} \label{dde.5}
\aligned
 &  \mathcal{D}(\mathcal{A}) = \{ x \in X \, : \, \dot{x} \in X \; \mbox{ and } \; \dot{x}(0) = A x(0) + B x(-\tau) \} \\
 &  \mathcal{A} \, x =  \dot{x}, \quad x \in \mathcal{D}(\mathcal{A}). 
\endaligned
\end{equation}
In this way, eq. (\ref{dde.3}) can be restated as the linear abstract Cauchy problem
\begin{equation} \label{dde.6}
\begin{aligned}
 &  \dot{U}(t) = \mathcal{A} \, U(t), \quad t \ge 0 \\
 & U(0) = \phi
\end{aligned} 
\end{equation}
on the Banach space $X$, where $U: [0, \infty) \longrightarrow \mathcal{D}(\mathcal{A})$ and the function $U(t) = x_t$ is the unique solution \cite{bellen03nmf}.

Typically, problem (\ref{dde.3}) (or (\ref{dde.6})) cannot be solved analytically, so that one has to find a way to get approximations. One possible way 
consists in discretizing the operator $\mathcal{A}$. This can be done by introducing a mesh
$\{ \theta_0, \theta_1, \ldots, \theta_N \}$ on the interval $[-\tau, 0]$, with $0=\theta_0 > \theta_1 > \cdots > \theta_N = -\tau$ and $N \ge 1$. Then, the
finite dimensional space $X_N = \mathbb{R}^{d(N+1)}$ is taken as the discretization of $X$. An element $\Phi = (\Phi_0, \ldots, \Phi_N)^T \in X_N$
can be seen as a function defined on the mesh, with $\Phi_j \in \mathbb{R}^d$, $j=0,1,\ldots, N$, being the value at $\theta_j$. As a simple example,
let us consider the equispaced mesh $\{ \theta_0 = 0, \theta_1=-h, \ldots, \theta_N=-N h \}$, with $h=\tau/N$. Then, the $d(N+1)$-dimensional
linear system of ODEs  
\begin{equation} \label{dde.7}
\begin{aligned}
 &  \dot{U}_N(t) = \mathcal{A}_N \, U_N(t), \quad t \ge 0 \\
 & U_N(0) = \phi_N
\end{aligned} 
\end{equation}
with $\phi_N = (\phi(0), \phi(-h), \ldots, \phi(-\tau))^T \in X_N$ and
\[
  \mathcal{A}_N = \left( \begin{array}{ccccc}
  			A & 0_d & \ldots & 0_d & B \\
			    &       & D_N \otimes I_d &  
			  \end{array} \right)  
\]
approximates the original system (\ref{dde.6}). Here $D_N$ is a $N \times (N+1)$ matrix corresponding to the particular finite difference scheme one chooses
for the integration, $I_d$ is the $d \times d$ identity matrix and $\otimes$ denotes the tensor product \cite{bellen00nso,bellen03nmf,maset99asi}. For instance,
if a first-order forward difference approximation is used, then
\[
   D_N = \frac{1}{h} \left( \begin{array}{cccccc}
                -1 & 1    &             &   &  & \\
                     &  -1 & 1          &   &  & \\
                     &       & \ddots & \ddots &    &  \\
                     &        &        & -1        & 1 &  \\
                      &       &        &   & -1 & 1 
                       \end{array} \right).        
\]                
A much more accurate description can be achieved by considering instead a pseudospectral differentiation method based on Chebyshev collocation points.
In this approach one takes the $N+1$ Chebyshev points
\[
  t_j =  \cos\frac{j\pi}{N}, \qquad j=0, \ldots, N
\]
on the interval $[-1,1]$ and the corresponding shifted points $\theta_j = (t_j -1) \tau/2$ as the mesh in $[-\tau, 0]$. Then, 
\begin{equation} \label{init.con}
   \phi_N = (\phi(\theta_0), \phi(\theta_1), \ldots, \phi(\theta_N))^T \in \mathbb{R}^{d(N+1)}
\end{equation}   
 and the matrix $\mathcal{A}_N$ is obtained as follows.
 First, one considers the standard 
$(N+1) \times (N+1)$ spectral differentiation matrix $D$ for the Chebyshev collocation points. The entries of $D$ can be found, e.g. in
\cite[p. 53]{trefethen00smi}. Then one forms the matrix $\mathbb{D} = D \otimes I_d$,  and
finally $\mathcal{A}_N$ is formed from $\mathbb{D}$ by replacing its first $d$ rows by zeros and replacing the
 $d \times d$ left upper corner by $A$
and the $d \times d$ right upper corner by $B$. In other words,
\begin{equation} \label{ma.A}
  \mathcal{A}_N = \left( \begin{array}{ccccc}
  			A & 0_d & \ldots & 0_d & B \\
			    &       & \frac{2}{\tau} [ \mathbb{D}^{(d+1,d(N+1))}] &  
			  \end{array} \right).  
\end{equation}
Here $[ \mathbb{D}^{(d+1, d(N+1))}]$ denotes the submatrix obtained by taking the rows $d+1, \ldots, d(N+1)$ from $\mathbb{D}$, whereas the
factor $2/\tau$ accounts for rescaling from the interval $[-1,1]$ to $[-\tau, 0]$ \cite{butcher11otc}. In particular, if $d=1$ and $N=4$, one has
\[
  \mathcal{A}_4= \frac{2}{\tau} \left(\begin{array}{ccccc}
       \frac{\tau}{2} A & 0 & 0 & 0 & \frac{\tau}{2} B \\
     1+ \frac{\sqrt{2}}{2}  &  -\frac{\sqrt{2}}{2} & -\sqrt{2} & \frac{\sqrt{2}}{2} & -\frac{1}{2 + \sqrt{2}} \\
     -\frac{1}{2} & \sqrt{2} & 0 & -\sqrt{2} & \frac{1}{2} \\
     \frac{1}{2+\sqrt{2}} & -\frac{\sqrt{2}}{2} & \sqrt{2} & \frac{\sqrt{2}}{2} & -1-\frac{\sqrt{2}}{2} \\
     -\frac{1}{2}  & \frac{4}{2+ \sqrt{2}} & -2 & \frac{4}{2- \sqrt{2}} & -\frac{11}{2}
     \end{array} \right)
\]     
where now $A$ and $B$ are scalars.
In this way, one ends up with the $d(N+1)$-dimensional initial value
problem defined by (\ref{dde.7}), with coefficient matrix (\ref{ma.A}) and initial condition (\ref{init.con}). By solving this system in the interval
$[0, \tau]$, one constructs a vector $U_N(\tau)$ whose entries constitute approximations to the solution of the DDE (\ref{dde.3}) at times 
$\theta_j^1 \equiv \tau + \theta_j$, namely 
$(U_N(\tau))_j \approx x(\tau + \theta_j)$. Thus, the first $d$ components of $U_N(\tau)$ provide the solution at $t = \tau$, the next $d$ components
approximate the solution at $t = \tau + \theta_1 = (1 + t_1) \frac{\tau}{2}$, etc.

If $N$ is sufficiently large, by applying this procedure one gets in fact spectral accuracy. This is observed, in particular,
when computing the eigenvalues of the approximate matrix $\mathcal{A}_N$: as shown in \cite{breda15sol}, the eigenvalues of $\mathcal{A}_N$ 
converge to the eigenvalues of the operator $\mathcal{A}$
faster than $\mathcal{O}(N^{-r})$ for any $r > 0$. Of course, the number of elements in the spectrum of 
$\mathcal{A}$ that are approximated by elements of the spectrum of $\mathcal{A}_N$ increases with the value of $N$, and the eigenvalues which are
closest of the origin are better approximated. Of course, the solution itself can be obtained by evaluating $U_N(\tau) = \exp(\tau \mathcal{A}_N) \phi_N$.

If one is interested in obtaining approximate solutions of the DDE (\ref{dde.3}) for $ t > \tau$, then the procedure consists in taking successive intervals
$[\tau, 2 \tau]$, $[2 \tau, 3 \tau]$, etc. and proceed as in the method of steps, taking as initial condition the approximation obtained in the previous interval.
Thus, if we denote $\theta_j^i \equiv  i \tau + \theta_j =  i \tau + \frac{(t_j-1)\tau}{2}$, to get approximations in the $(i+1)$-th interval $[i \tau, (i+1) \tau]$,
$i \ge 1$, we have to solve (\ref{dde.7}) with the initial condition $U_N(i \tau)$ obtained by integrating on the previous $i$-th interval. If, on the other hand,
we want to get approximations at a particular time $t \in ( i \tau, (i+1) \tau)$, then an interpolation can be carried out based on the values obtained at the
Chebyshev points in the interval.

\subsection{Linear periodic DDEs}

Whereas in the autonomous case, the infinitesimal generator approach allows one to get accurate approximations to the eigenvalues and the solution
of (\ref{dde.1}), and even obtain rigorous convergence estimates, the situation is far more complicated 
when matrices $A$ and $B$ in (\ref{dde.1}) depend explicitly on time. One could analogously try to describe the time evolution of the linear DDE through
a Cauchy problem defined for an abstract ODE, but in that case 
one has to introduce a 2-parameter family of operators called \emph{evolution family} as solutions for abstract Cauchy problems of the form
\cite[Chapter 3]{chicone99esi}
\begin{equation} \label{lp.1}
  \dot{y}(t) = \mathcal{A}(t) \, y(t), \quad y(s) = z \in \mathcal{D}(\mathcal{A}(s)), \quad t \ge s,
\end{equation}
where the domain $\mathcal{D}(\mathcal{A}(s))$  of the operator $\mathcal{A}(s)$ is assumed to be dense in $X$. Specifically, a family 
$\{U(t,s)\}_{t \ge s}$ of linear and bounded operators on $X$ is called an evolution family if
\begin{itemize}
 \item[(i)] $U(t,s) = U(t, \nu) U(\nu,s)$ and $U(s,s) = I$ for all $t \ge \nu \ge s$; and
 \item[(ii)] for each $\phi \in X$, the function $(t,s) \mapsto U(t,s) \phi$ is continuous for $t \ge s$.
\end{itemize}
The abstract Cauchy problem (\ref{lp.1}) is called \emph{well-posed} if there exists an evolution family $\{U(t,s)\}_{t \ge s}$ that solves (\ref{lp.1}), i.e., if
for each $s \in \mathbb{R}$, there exists a dense subset $Y_s \subseteq \mathcal{D}(\mathcal{A}(s))$ such that, for each $y_s \in Y_s$ the function
$t \mapsto y(t) \equiv U(t,s) y_s$, for $t \ge s$, is differentiable, $y(t) \in \mathcal{D}(\mathcal{A}(t))$ and (\ref{lp.1}) holds \cite{chicone99esi,engel00ops}.

It has been shown that the non-autonomous Cauchy problem (\ref{lp.1}) is well-posed if and only if there exists a unique evolution family
$\{U(t,s)\}_{t \ge s}$  solving (\ref{lp.1}) \cite{engel00ops,nickel97esf}. If the well-posedness of the problem has been established (which in some cases
is far from trivial \cite{chicone99esi,engel00ops,csomos22aso}), then the unique solution has the form
\[
   y(t) = U(t,s) z \qquad \mbox{ for } \qquad t \ge s.
 \]


In the particular case of periodic problem with period $T$,
\begin{equation} \label{lp.2}
\aligned
  & \dot{x} = L(t) \,  x_t, \quad t \ge s \\
  & x_s = \phi,
\endaligned  
\end{equation}
the evolution operator $U(t,s) \phi = x_t(\cdot; s, \phi)$ verifies in addition that
$U(t+T,s) = U(t,s) U(s+T,s)$ for all $t \ge s$, and it is the so-called monodromy operator $U(T,0)$ what is the central
object of study to determine the stability of the system. It can be shown that the spectrum of $U(T,0)$ is an at most countable compact set of $\mathbb{C}$
with zero as the only possible accumulation point. Moreover, any eigenvalue $\mu \ne 0$ belongs to the point spectrum, and is called a characteristic
multiplier of eq. (\ref{dde.1}) \cite{hale93itf}. In particular, the zero solution of (\ref{lp.2}) is uniformly asymptotically stable if and only if all the characteristic
multipliers are such that $|\mu| < 1$ \cite{breda15sol,hale93itf}. 
In any event, the solution $y(t)$ (corresponding to $x(t)$ in (\ref{dde.1})) has to be approximated at certain values of $t$.

Whereas in \cite{breda15sol} a pseudospectral collocation method is applied directly to discretize the monodromy operator $U(T,0)$ and to compute
approximations to the characteristic multipliers of (\ref{dde.1}), in this paper we proceed  formally as in the case of autonomous problems. In other
words, we take $N+1$ Chebyshev nodes in the interval $[-\tau, 0]$, replace $X$ by the finite dimensional space $X_N$ and the abstract Cauchy problem 
by
\begin{equation} \label{lp.3}
\begin{aligned}
 &  \dot{U}_N(t) = \mathcal{A}_N(t) \, U_N(t), \quad t \ge 0 \\
 & U_N(0) = \phi_N,
\end{aligned} 
\end{equation}
where now 
\begin{equation} \label{lp.4}
  \mathcal{A}_N(t) = \left( \begin{array}{ccccc}
  			A(t) & 0_d & \ldots & 0_d & B(t) \\
			    &       & \frac{2}{\tau} [ \mathbb{D}^{(d+1,d(N+1))}] &  
			  \end{array} \right)  
\end{equation}
is a periodic matrix of period $T$ and the initial vector is the discretization of the function $\phi(t)$ at the shifted Chebyshev points,
\begin{equation} \label{lp.5}
  \phi_N = (\phi(\theta_0), \phi(\theta_1), \ldots, \phi(\theta_N))^T \in \mathbb{R}^{d(N+1)}. 
\end{equation}
  This problem is then numerically integrated to get approximations to the solution over the interval $[0,\tau]$, specifically at the times determined by the Chebyshev
points. Depending on the value of $N$ we get approximations at more points in $[0,\tau]$. Once the
approximation at $t=T$ is obtained, then the monodromy operator is available, so that one readily determines the first $d(N+1)$ characteristic
multipliers of eq. (\ref{dde.1}) \cite{butcher11otc}. 

\section{Numerical integration of linear DDEs}
\label{sec.3}

\subsection{Numerical integrators based on the Magnus expansion}
\label{sec.3.1}

One is then confronted with the numerical time integration of the non-autonomous 
initial value problem (\ref{lp.3})-(\ref{lp.5}) defined in the finite dimensional space $X_N$.
This of course can be done by applying
several numerical integrators, such as Runge--Kutta or multistep methods. Among them, numerical methods based on the Magnus expansion are
particularly appropriate for linear systems, resulting in very efficient schemes that, in addition, preserve important qualitative properties of the continuous system
\cite{iserles00lgm}. It makes sense, then, trying to combine the spectral accuracy provided by the Chebyshev collocation points with Magnus integrators to get
accurate approximations.

Magnus' approach to solve the general linear differential equation
\begin{equation} \label{me.1}
  \dot{y}(t) = \hat{A}(t) \, y(t), \qquad y(0) = y_0,
\end{equation}  
where $\hat{A}(t)$ is a $m  \times m $ matrix and $y \in \mathbb{R}^m$, consists in expressing the solution as an exponential $y(t) = \exp(\Omega(t)) y_0$,
and determining the equation satisfied by $\Omega(t)$. Specifically, it can be shown that
\begin{equation} \label{me.2}
  \dot{\Omega} = \sum_{k=0}^{\infty} \frac{B_k}{k!} \ad_{\Omega}^k (\hat{A}(t)), \qquad \Omega(0) = 0,
\end{equation}
where $B_k$ are the Bernoulli numbers and $\ad_{\Omega}(\hat{A}) = [ \Omega, \hat{A}] = \Omega \, \hat{A} - \hat{A} \, \Omega$, 
$\ad_{\Omega}^k (\hat{A}) = [\hat{A}, \ad_{\Omega}^{k-1} (\hat{A})]$.  Equation (\ref{me.2})
is then solved  by applying Picard fixed point iteration after integration. This results in an infinite series for $\Omega$, 
 \begin{equation}
\Omega(t)=\sum_{k=1}^{\infty}\Omega_{k}(t), \qquad \mbox{ with } \qquad \Omega_k(0) = 0,  \label{me.3}%
\end{equation}
whose terms are increasingly complex expressions involving time-ordered integrals  of nested commutators of $\hat{A}$ evaluated at different times
\cite{blanes09tme,arnal18agf}. In particular,
\begin{align} \label{ma.4}
\Omega_{1}(t)  & =\int_{0}^{t}\hat{A}(t_{1}) dt_{1},\nonumber\\
  \Omega_{2}(t) &  = - \frac{1}{2}\int_{0}^{t} dt_{1} \int_{0}^{t_{1}}
dt_{2}\ \left[  \hat{A}(t_{2}), \hat{A}(t_{1})\right]. 
\end{align}
The expansion is guaranteed to converge at least for $t \in [0, t_c)$ such that
\[
  \int_0^{t_c} \|\hat{A}(s)\|_2 \, ds < \pi.
\]  
By appropriately truncating the series and approximating the integrals by suitable quadratures, it is then possible to design numerical integration schemes 
for approximating the solution of eq. (\ref{me.1}) in a given interval $[0, t_f]$ \cite{blanes09tme,iserles00lgm}. 
As usual, the integration interval is divided into $M$ steps such that (\ref{me.3}) converges in
each subinterval $[t_{k}, t_{k+1}]$, $k=0,\ldots M$, with $t_M = t_f$ and length $h = t_f/M$. 
Then, an approximation to the solution of (\ref{me.1}) of order $2p$ is achieved by computing
\begin{equation} \label{num.int1}
  y(t_{k+1}) \approx  y_{k+1} = \exp( \Omega^{[2p]}(h)) \, y_{k},  \qquad k=0,1, \ldots, M-1,
\end{equation}
where $\Omega^{[2p]}(h) = \Omega(h) + \mathcal{O}(h^{2p+1})$. Since the construction procedure
is detailed in references \cite{blanes00iho,blanes09tme}, here we only collect specific integration schemes of order $2p= 2, 4, 6$.

\paragraph{Order 2.} 
A 2nd-order scheme is attained by approximating $\Omega_1$ by the midpoint rule. It reads
\begin{equation} \label{or2}
   y_{k+1} = \exp( h \hat{A}(t_k + h/2)) \, y_k, \qquad k=0,\ldots, M-1
\end{equation}      
and is also known as the exponential midpoint rule.

\paragraph{Order 4.}
If we take the 2 points Gauss--Legendre quadrature rule and compute
\[
  A_1= \hat{A} \left( t_k + \left( \frac{1}{2} - \frac{\sqrt{3}}{6} \right) h \right), \qquad   A_2 = \hat{A} \left( t_k + \left( \frac{1}{2} + \frac{\sqrt{3}}{6} \right) h \right),
\]
then the scheme
\begin{equation} \label{or4}
\aligned
  & \Omega^{[4]}(h) = \frac{h}{2} (A_1 + A_2) - h^2 \frac{\sqrt{3}}{12} \, [A_1, A_2] \\
  & y_{k+1} = \exp( \Omega^{[4]}(h)) \, y_k
\endaligned
\end{equation}  
renders an integration method of order 4 for equation (\ref{me.1}).

\paragraph{Order 6.}
By evaluating the matrix $\hat{A}$ at the Gauss--Legendre collocation points
\[
  A_1= \hat{A}(t_k + c_1 h), \qquad A_2 = \hat{A}(t_k + c_2 h), \qquad A_3 = \hat{A} (t_k + c_3 h),
\]
with $c_1 = 1/2 - \sqrt{15}/10$, $c_2 = 1/2$,   $c_3 = 1/2 + \sqrt{15}/10$, we form the quantities
\[
  \alpha_1 = h A_2, \qquad \alpha_2 = \frac{\sqrt{15} h}{3} (A_3 - A_1), \qquad \alpha_3 = \frac{10 h}{3} (A_3 - 2 A_2 + A_1),
\]
and finally
\begin{equation} \label{or6}
\aligned
  & C_1 = [\alpha_1, \alpha_2] \\
  & C_2 = -\frac{1}{60} [\alpha_1,  2 \alpha_3 + C_1] \\
  & \Omega^{[6]}(h) = \alpha_1 + \frac{1}{12} \alpha_3 + \frac{1}{240} [-20 \alpha_1 - \alpha_3 + C_1, \alpha_2 + C_2] \\
  & y_{k+1} = \exp( \Omega^{[6]}(h)) \, y_k.
\endaligned
\end{equation}  
Then we end up with an integrator of order 6 requiring only the evaluation of $\hat{A}$ at three different times and the computation of 3 commutators.

These schemes have shown to be more efficient than standard integrators such as Runge--Kutta methods, even when the evaluation
of the exponential of a matrix is required at each step \cite{blanes00iho}. In this respect, we note that several  algorithms have been recently
proposed to reduce the computational cost of the exponential matrix \cite{bader19ctm}, especially in the context of exponential integrators \cite{bader22aea}.
These will be incorporated into our procedure.

\subsection{Integration algorithm}
\label{sec.3.2}

We have now all the tools required to formulate a practical algorithm for the numerical integration of the linear DDE (\ref{dde.1}) in the interval 
$t \in [0, \tau]$, for generic
(sufficiently smooth) $d \times d$ matrices $A(t)$, $B(t)$ and an initial function $\phi(t)$. It can be summarized as follows:
\begin{enumerate}
 \item Choose an integer $N$ and form the spectral differentiation matrix $\mathbb{D}$ associated with the $N+1$ Chebyshev collocation points on the
 interval $[-1,1]$. 
 \item Construct the $d(N+1) \times d(N+1)$ matrix $\mathcal{A}_N(t)$ of eq. (\ref{lp.4}).
 \item Form the $d(N+1)$-dimensional vector $U_N(0)$, whose elements are the values of the initial function $\phi$ at the shifted Chebyshev collocation 
 points on the interval $[-\tau,0]$.
 \item Choose an integer $M$, the step size $h = \tau/M$ and construct the numerical solution $\widetilde{U}_N(M h)$ of the initial value problem 
 (\ref{lp.3})-(\ref{lp.5}) at the final time
$t = \tau$ by applying one
of the numerical integrators of section \ref{sec.3.1} obtained from the Magnus expansion.
\item The vector $\widetilde{U}_N(M h)$ provides an approximation to $U_N(\tau)$ of order of convergence $2p$, i.e., 
\[
    \widetilde{U}_N(M h) = U_N(\tau) + \mathcal{O}(h^{2p+1}), \qquad p=1,2,3,
\]
and thus also approximates the solution of eq. (\ref{dde.1}) at times $\tau + \theta_j$, $j=0, \ldots, N$.
\end{enumerate}

As pointed out before, if the solution has to be computed at times $t > \tau$, the same process can be applied in any interval $[i \tau, (i+1) \tau]$, $i \ge 1$,
taking as initial condition the approximation at the end of the previous interval, $\widetilde{U}_N(M h)$. The matrix $\mathcal{A}_N(t)$ 
has to be formed again, of course, but the only difference with respect to the previous interval lies in the time-dependent part
 (i.e., in $A(t)$ and $B(t)$).
 
 The procedure has thus three main ingredients: (a) the original DDE is reformulated as an abstract Cauchy problem in the Banach space
 $X$; (b) this problem is then discretized by a pseudo-spectral method, thus leading to an initial value problem defined in a finite-dimensional space $X_N$,
 and (c) finally a Magnus integrator is applied to solve numerically the resulting finite-dimensional non-autonomous ordinary differential equation.
 In this sense, the algorithm thus combines the advantages of both spectral and Magnus methods.
 
 Notice, in particular, that the procedure critically depends on the parameters $N$ (related to the discretization of the delayed part) and $M$ (leading to the step size 
 $h$ in
 the time numerical integration), and finally also on the order of accuracy of the Magnus integrator. These parameters can
 be chosen according with the required accuracy, and are clearly related. The situation closely resembles what happens when a Magnus integrator is used
 to integrate in time a partial differential equation previously discretized in space: one has to adjust $M$ to provide an accuracy consistent with
the scheme used for the space discretization.

In this respect it is worth noticing that the Magnus expansion provides the exact solution of the initial value problem (\ref{lp.3})-(\ref{lp.5}) in the autonomous
case (i.e., when the matrices $A$ and $B$ in (\ref{dde.1}) are constant). This feature can be used to choose the number of collocation points $N$ leading to
the required accuracy in the general case simply by freezing $t$ at one particular time. Once a particular $N$ has been selected, the number of subdivisions
$M$ and the order of the Magnus integrators can be fixed according with some specified tolerance. The order and the step size can even be changed 
from one particular interval $[ i \tau, (i+1) \tau ]$ to the next \cite{blanes02hoo}. 

If the number of nodes $N$ of the spectral discretization is
sufficiently large, then the algorithm has the usual properties exhibited by standard Magnus integrators applied to non-autonomous linear differential equations
concerning stability, convergence and error propagation \cite{blanes09tme,blanes00iho,iserles00lgm}.

 In the periodic case, one can use this procedure to compute the fundamental matrix of system (\ref{lp.3}) at time $t=T \ge \tau$, and its eigenvalues. Then, they 
provide approximations to the first $d(N+1)$ characteristic multipliers of the DDE (\ref{dde.1}). As in the autonomous case, the eigenvalues that are closest to
its minimum value (in absolute value) are better approximated.

It is worth remarking that, whereas the previous algorithm only deals with one delay, it can be easily generalized to the case of multiple discrete delays, i.e.,
to solve numerically equation (\ref{dde.1b}) by applying any of the alternatives proposed by \cite{butcher11otc}: either by fitting a single Chebyshev
polynomial through the entire delay interval $[- \tau_k, 0]$ or by considering a different Chebyshev polynomial of different degree in every interval 
$[-\tau_k, -\tau_{k-1}]$, $[-\tau_{k-1}, -\tau_{k-2}]$, etc.

\section{Numerical examples}
\label{sec.4}

Next we illustrate the previous algorithm in practice on several examples, both for computing 
the first characteristic multipliers of the problem and also for getting 
numerical approximations to the exact solution. For convenience, we denote the previous Magnus
integrators of order 2, 4 and 6 as M2, M4 and M6, respectively.

\paragraph{Example 1: a scalar periodic equation.}
As a first illustration we take the admittedly simple non-autonomous one-dimensional equation with $2 \pi$-periodic coefficients and delay $\tau = \frac{\pi}{2}$
\begin{equation} \label{1d.1}
 \dot{x}(t) = \cos(t) \,x(t) - \e^{\sin t + \cos t} \, x(t-\frac{\pi}{2}).
\end{equation}
It is taken from \cite{breda15sol} and arises when linearizing the nonlinear autonomous DDE
\begin{equation} \label{1d.2}
  \dot{z}(t) = -\log \left( z\left(t-\frac{\pi}{2}\right) \right) \, z(t).
\end{equation} 
around its periodic solution $z(t) = \exp(\sin t)$. It is a simple exercise to check that $x(t) = \dot{z}(t) = \e^{\sin t} \cos t$ is indeed a solution of
(\ref{1d.1}), so that $\mu = 1$ is a characteristic multiplier \cite{breda15sol}. We can therefore check the accuracy of the algorithm by computing both
the spectrum of the monodromy matrix and the numerical approximation to the exact solution by integrating eq. (\ref{1d.1}) with initial function
$\phi(t) = \e^{\sin t} \cos t$.

\begin{figure}[!ht] 
\begin{center}
\includegraphics[width=6.5cm]{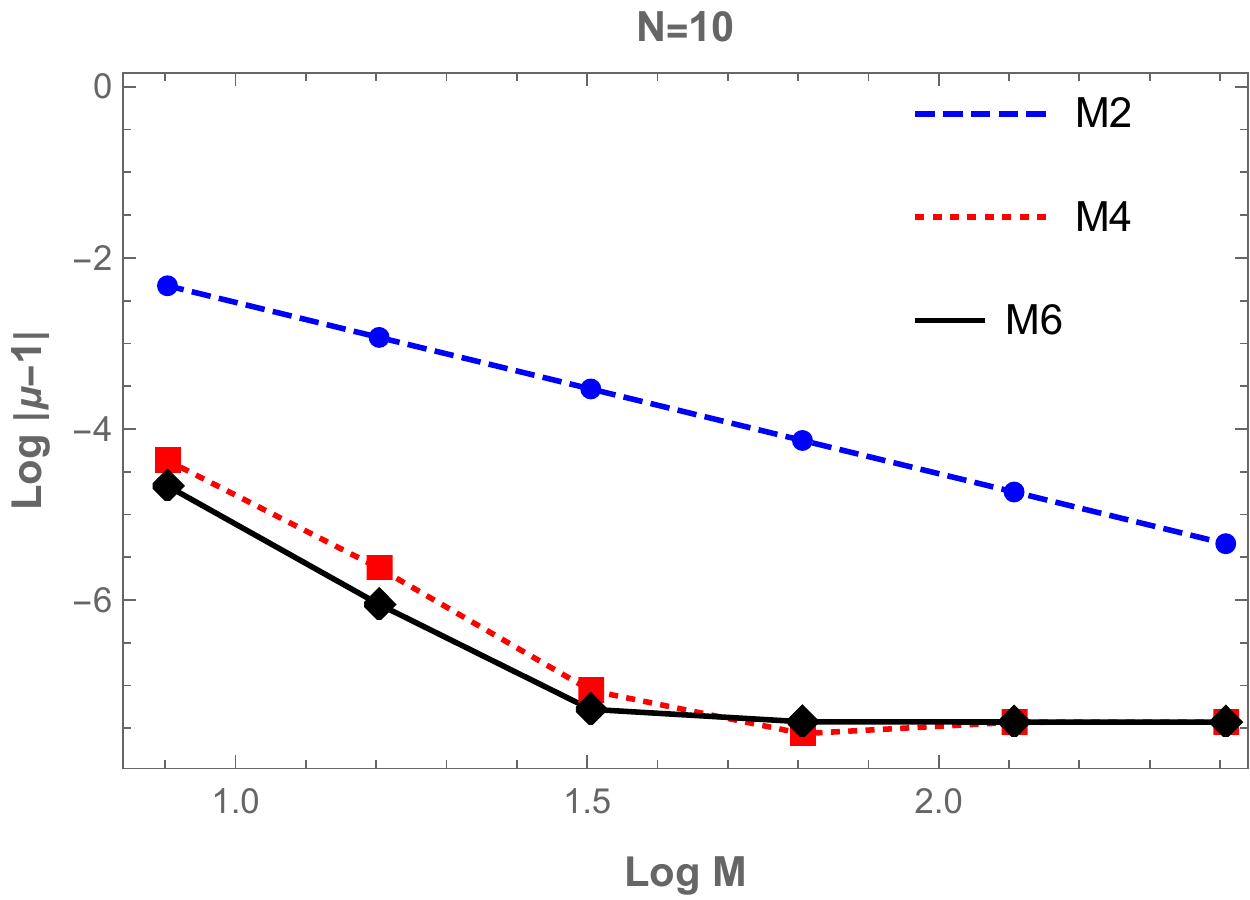} 
\includegraphics[width=6.6cm]{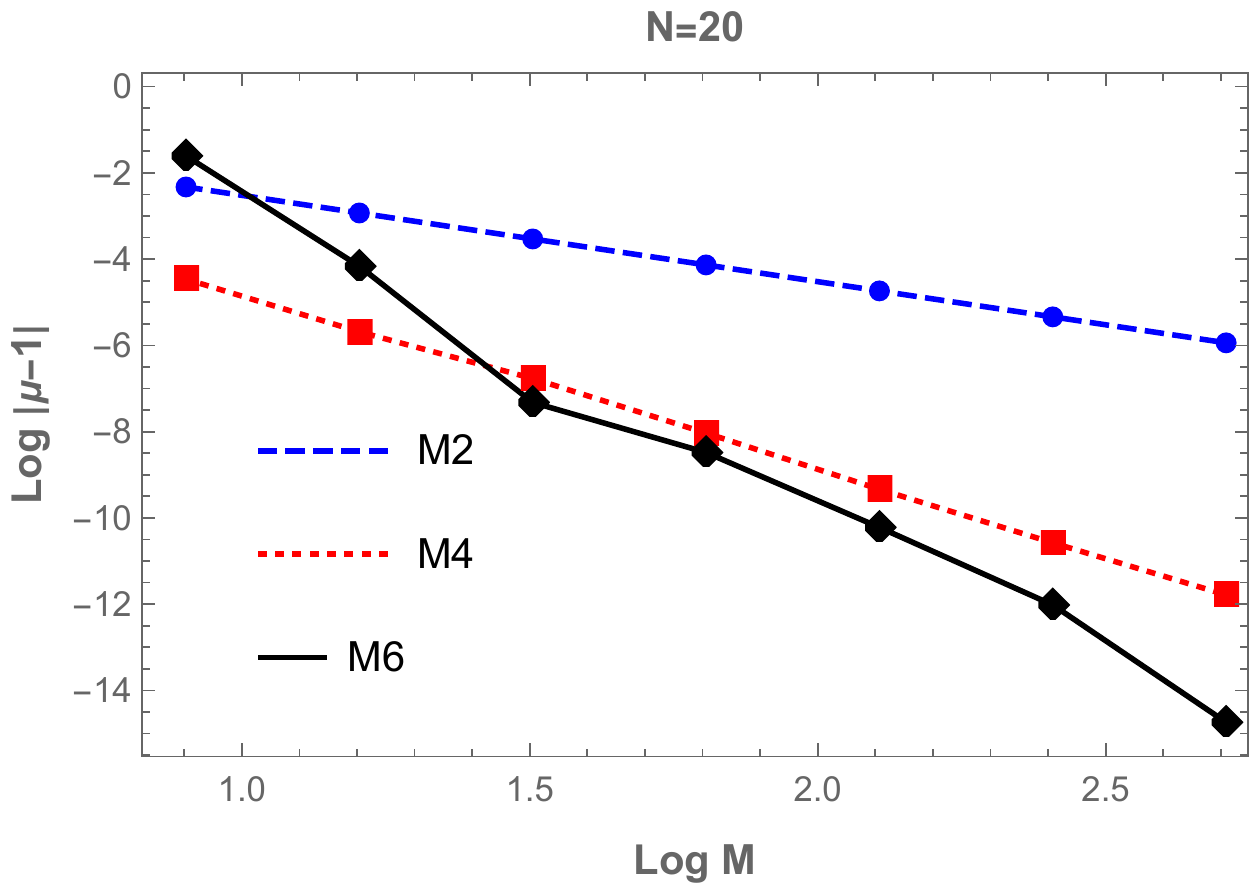} 
\end{center}
\caption{\label{figure1s} \small Error in the dominant characteristic multiplier $\mu = 1$ for the problem (\ref{1d.1}) obtained with Magnus integrators M2, M4
and M6  as a function of $M$ with two different values for the number of collocation points: $N=10$ (left) and $N=20$ (right).}
\end{figure}

The first test is devoted to check the accuracy in the determination of the characteristic multiplier $\mu = 1$. This is done by computing the monodromy matrix
through the numerical integration of the matrix system 
\[
  \dot{Y} = \mathcal{A}_N(t) Y, \qquad Y(0) = I
\]
associated with (\ref{lp.3})-(\ref{lp.4})  until the final time $t_f = 2 \pi$, i.e., in the interval $[0, 4 \tau]$, with the previous Magnus integrators, and then
by  obtaining its eigenvalues . We then compute the difference between the first eigenvalue and $\mu = 1$ as a function of $M$ (the number 
of subdivisions in each interval $[i \tau, (i+1)\tau]$) for different values of $N$ (the number of collocation points). In this way, we obtain the results collected
in Figure \ref{figure1s}, with $N=10$ (left diagram) and $N=20$ (right panel). 
The order of approximation of M2 and M4 is clearly visible in the figure, whereas the higher accuracy of M6 is only visible for sufficiently large $N$ and $M$.

Notice that the main limiting factor for accuracy in the left panel is the small number of collocation
points: doubling the value of $N$, from 10 to 20, allows us to decrease the error by more than 5 orders of magnitude. In this respect, it is worth 
mentioning that with the procedure proposed in \cite{breda15sol}, $N=40$ collocation points are required to render similar errors as those achieved here
with $N=20$.

\begin{figure}[!ht] 
\begin{center}
\includegraphics[width=6.5cm]{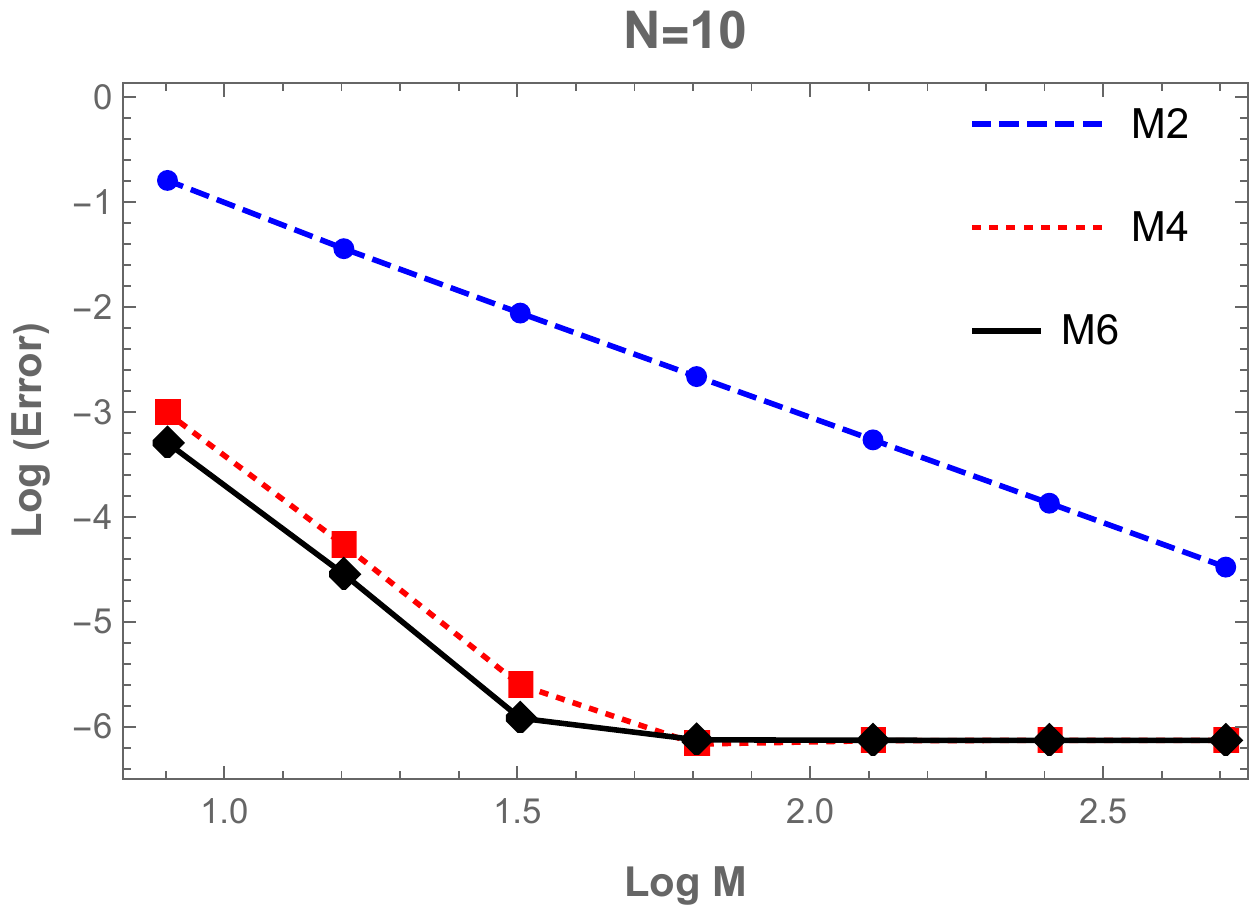} 
\includegraphics[width=6.6cm]{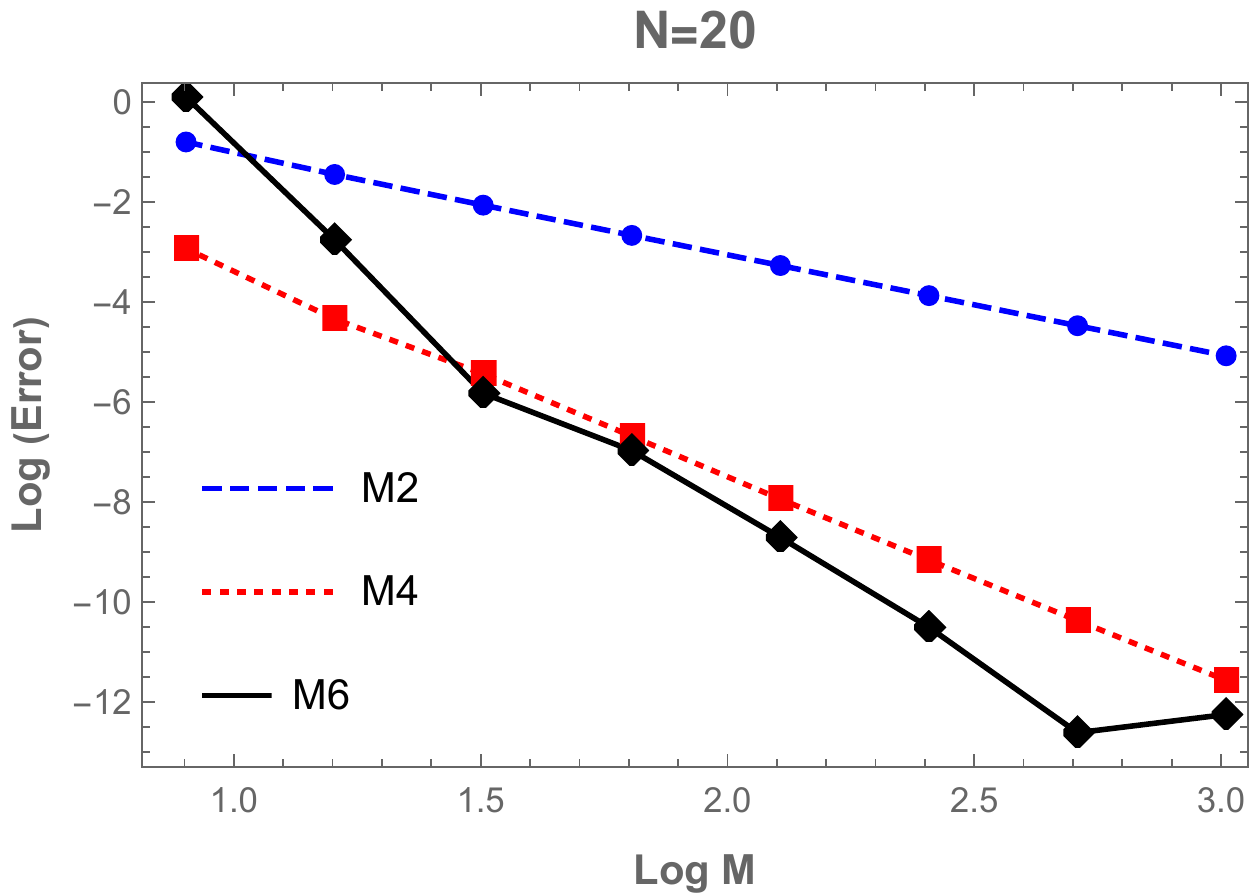} 
\end{center}
\caption{\label{figure2s} \small Mean error (\ref{mean-error}) in the solution $x(t) = \e^{\sin t} \cos t$ of (\ref{1d.1}) obtained with Magnus integrators M2, M4
and M6  as a function of $M$ with two different values for the number of collocation points: $N=10$ (left) and $N=20$ (right).
Final time: $t_f = 100 \pi$.}
\end{figure}

In our second experiment we check the accuracy in solving the problem (\ref{1d.1}) for long times by computing approximations to the exact solution.
Specifically, we integrate until the much larger final time $t_f = 100 \pi$ with M2, M4 and M6 for different values of $M$ and compute 
the mean error of the solution in
the last interval $[199 \tau, 200 \tau]$,
\begin{equation} \label{mean-error}
  \mathcal{E} = \frac{1}{N+1} \sum_{j=0}^N |x(t_j) - (U_N)_{j}|.
\end{equation}
Here $x(t_j) = \e^{\sin t_j} \cos t_j$, $t_j = 200 \tau + \theta_j$  and $U_N$ are the approximations obtained with each integrator. 
The results are displayed in Figure \ref{figure2s}
for $N=10$ (left) and $N=20$ (right) collocation points. The same notation as in Figure \ref{figure1s} has been used for each method. Here again, the higher
order of M6 is already visible with only $N=20$ collocation points, and the quality of the approximation does not degrade even if long time integrations
are considered. Notice the close similarity between the results exhibited in Figures \ref{figure1s} and \ref{figure2s}.

\paragraph{Example 2: delayed Mathieu equation.}
The DDE we consider next is
\begin{equation} \label{dme1}
  \ddot{x}(t) + ( \delta + \varepsilon \cos t) x(t) = b \, x(t - \tau), 
\end{equation}
where $\delta$, $\varepsilon$ and $b$ are real parameters.
The case in which the time delay $\tau$ is equal to the principal period $2 \pi$ has been thoroughly studied in the literature, especially with respect to 
its stability \cite{insperger11sdf}, and so we also fix $\tau = 2 \pi$ here. Equation (\ref{dme1}) includes both the effects of time delay and the presence of
parametric forcing, and appears in relevant mechanical engineering problems (see  \cite{insperger11sdf} and references therein).

\begin{figure}[!ht] 
\begin{center}
\includegraphics[width=6.5cm]{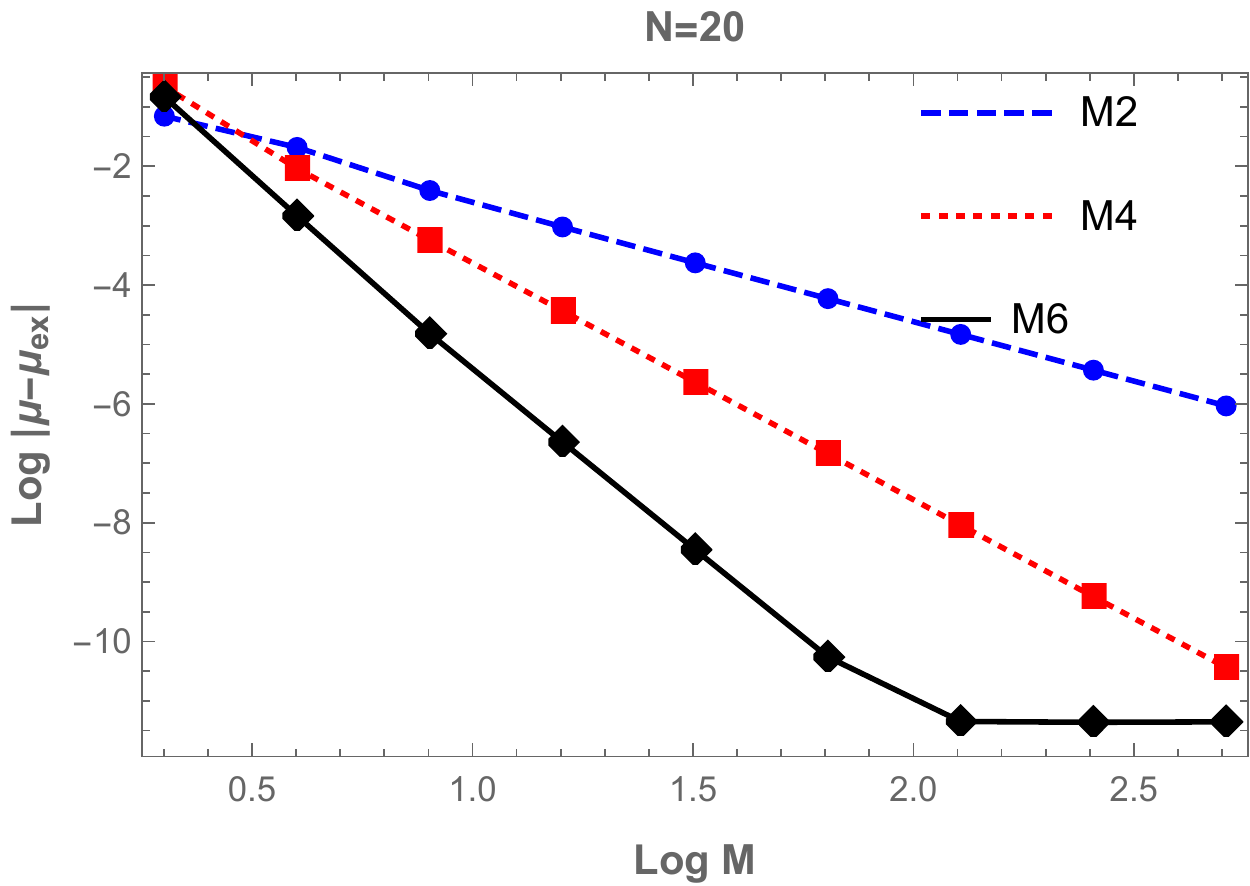} 
\includegraphics[width=6.5cm]{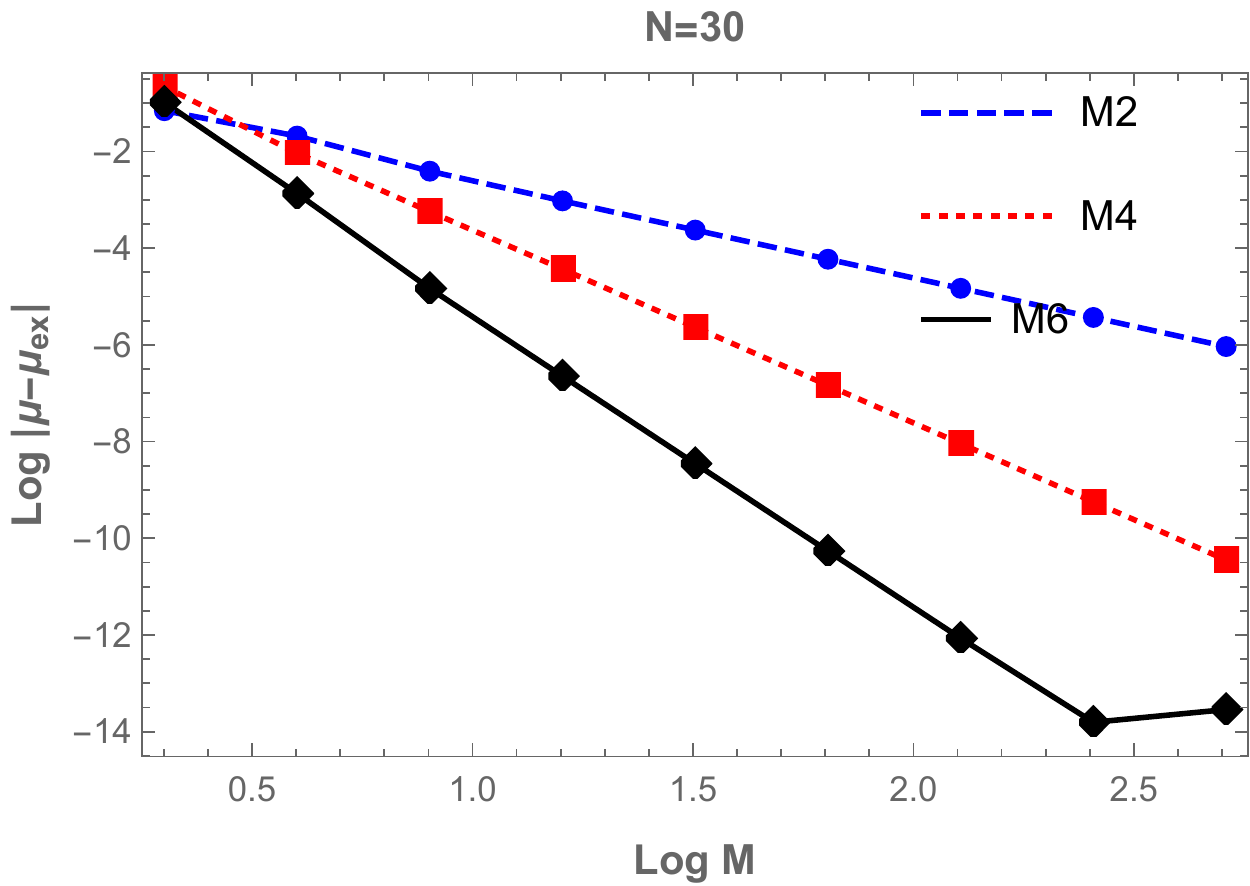} 
\end{center}
\caption{\label{figure3} \small Error in the first characteristic multiplier $\mu_{ex} $ for the delayed Mathieu equation (\ref{dme1})
obtained with Magnus integrators M2, M4
and M6 as a function of $M$ for $N=20$ (left) and $N=30$ (right). The values of the parameters are: $\delta = 1.5$, $\varepsilon = 0.5$, $b=-0.2$.}
\end{figure}

Equation (\ref{dme1}) is transformed into a system of the form (\ref{dde.1}) with $d=2$ and matrices
\[
   A(t) =  \left( \begin{array}{cr}
                    0 & \  1 \\
                    -(\delta + \varepsilon \cos t) & \ 0
                    \end{array}  \right),  \qquad
  B =   \left( \begin{array}{cc}
          b & 0 \\
          0 & 0
          \end{array} \right),
\]   
to which the previous algorithm can be readily applied. In our first test we fix the parameters to the values $\delta = 1.5$, $\varepsilon = 0.5$ and $b=-0.2$.
By using the Floquet technique proposed in \cite{insperger02scf}, it can be seen that one of the characteristic multipliers is given by
\[
  \mu_{ex} = 0.22751840350292177638239482513 + 1.417175174215530683457881875737 \, i
\]
with 30 digits of accuracy. This is taken as the reference value to compare with our procedure. As in the previous example, we compute the corresponding
monodromy matrix with Magnus integrators M2, M4 and M6 for different values of the discretization parameters $N$ and $M$. We then determine its
first eigenvalue $\mu_1$ and the error $|\mu_1 - \mu_{ex}|$ by fixing $N$ and several values of $M$. The corresponding results are depicted in 
Figure \ref{figure3} (in a log-log scale) when $N=20$ (left) and $N=30$ (right) as a function of $M$. Again, the order of each scheme is clearly visible,
with M6 providing the most accurate approximations for all the values of $M$ considered. Here also the number of collocation points $N$ determines the
maximal accuracy one can get. For comparison, with the numerical technique proposed in \cite{breda15sol} (based on the direct discretization of the
evolution operator) one gets also round-off accuracy with $N=30$,
whereas the error achieved with $N=20$ is $\approx 3.1 \cdot 10^{-9}$. Notice that, even with a time step as large as $h=\tau/2$  (or $M=2$), there are
no stability issues in the results.  

To check how the errors in the numerical solution evolve with time, for our next experiment we take as initial condition $\phi(t) = t$ and integrate 
the initial value problem (\ref{lp.3})-(\ref{lp.5}) until the final time $t_f = 10 \tau = 20 \pi$ with a fixed value for both discretization parameters, $N=40$,
$M=200$. As a measure of the error we compute the quantity
\[
  E_i = \frac{1}{N+1} \sum_{j=0}^{N-1} |x(t_{j}) - (U_N)_{2j+1}|, \qquad \mbox{ with } \qquad t_j \in [i \tau, (i+1) \tau]
\]
in each subinterval, where the exact solution is taken as the output of
the function \texttt{NDSolve} of \emph{Mathematica} with a very stringent tolerance.
 Notice that, since $d=2$, the odd components of the vector $U_N$ provide the approximations to the exact solution of (\ref{dme1}) at
the collocation points, whereas the elements $(U_N)_{2j}$ approximate $\dot{x}(t_j)$.  

\begin{figure}[!ht] 
\begin{center}
\includegraphics[scale=0.8]{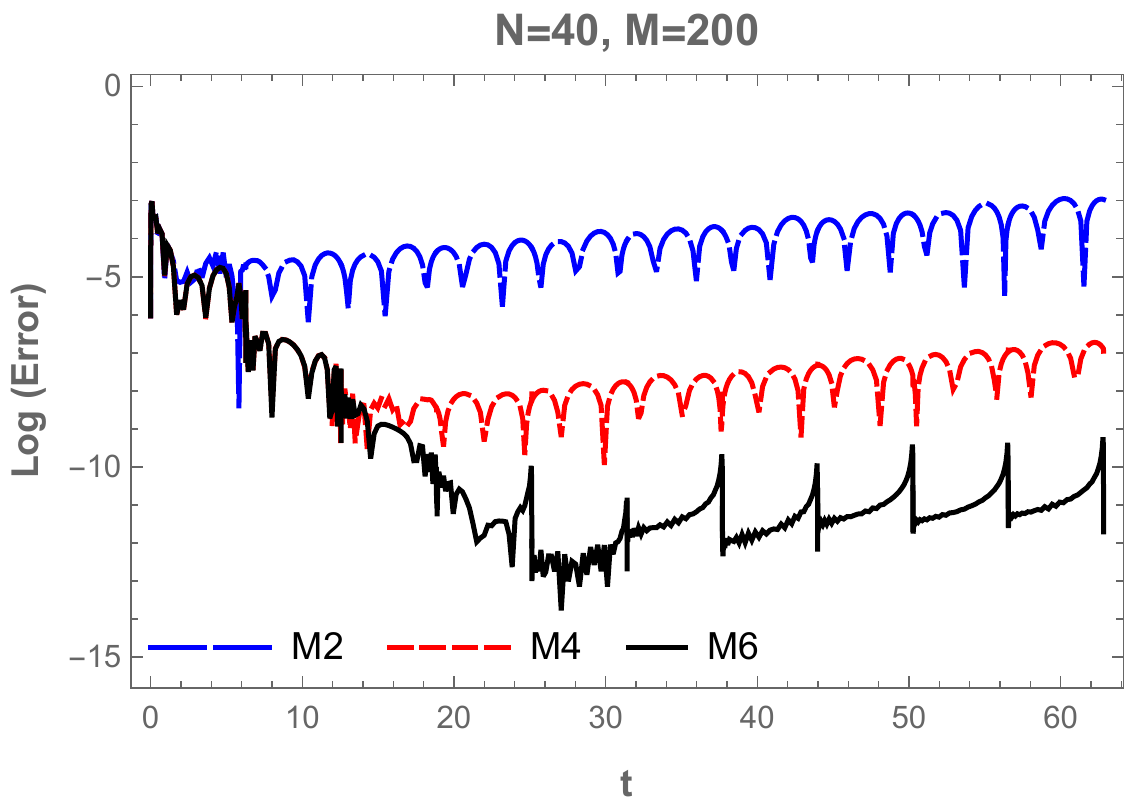} 
\end{center}
\caption{\label{figure4} \small Mean error in the solution of the delayed Mathieu equation (\ref{dme1}) with $\delta = 1.5$, $\varepsilon = 0.5$, $b=-0.2$
and initial function $\phi(t) = t$ obtained by M2, M4 and M6 in the interval $[0, 20 \pi]$ with $N=40$, $M=200$.}
\end{figure}

We observe that initially the errors achieved by M4 and M6 are similar to those committed by M2. This feature can be explained by the fact that 
 the initial function does not satisfy the differential equation. After
several subintervals, however, the accuracy actually improves and remains so for the integration interval considered.

As a final test, we illustrate how our numerical algorithm can also be used to determine with high accuracy the stability of a system modeled by equation
(\ref{dme1}), and in particular the stability boundaries in the parameter space, such as are presented in the stability chart of 
\cite[Fig. 2.10]{insperger11sdf}. To this end we fix $\delta = 2$ and $\varepsilon = 1$. Then it can be shown that there exists a value of $b \approx 0.71$ for which
the dominant multiplier is $\mu = 1$. Specifically, the spectral collocation method of \cite{breda15sol} with $N=20$ provides $b=0.706833720464083$.
That this is not the correct value (up to machine accuracy) can be seen either by increasing $N$ with the same algorithm or by applying the technique
of \cite{insperger02scf}, resulting instead in $b=0.7068337166604264$. For this value of $b$, the algorithm proposed here with $N=20$ and the 6th-order
Magnus integrator with $M=40$ provides the characteristic multiplier with an error $\approx 5.34 \cdot 10^{-12}$, in contrast with an error $\approx 5.38 \cdot
10^{-9}$ for the procedure of \cite{breda15sol} also with the same number of collocation points $N=20$.

\section{Extension to quasilinear problems}
\label{sec.5}

The algorithm we have presented in the previous section is mainly addressed to linear DDEs, both autonomous and explicitly time-dependent. 
It turns out, however, that the same procedure
can also be formally extended to the quasilinear delay equation (\ref{dde.2}), which we write again here for convenience:
\begin{equation} \label{eq.5.1}
\aligned
  & \frac{dx}{dt} = A \big(x(t-\tau) \big) \, x(t), \quad t \ge 0 \\
  & x(t) = \phi(t), \qquad -\tau \le t \le 0.
\endaligned
\end{equation}
This problem can be formally stated as 
\[
 \aligned
  & \frac{d}{dt} \tilde{u}(t) = Q(t) \tilde{u}(t), \qquad t \ge 0  \\
  & \tilde{u}(0) = \phi(0)
 \endaligned
\]
in terms of the function $\tilde{u}: [-\tau, \infty) \longrightarrow X$
\[
  \tilde{u}(t) \equiv \left\{ \begin{array}{lr}
  		\phi(t), & t \in [-\tau, 0] \\
		x(t), & t \in [0, \infty)
	\end{array} \right.
\]
and the operator $Q(t) \equiv A(x(t-\tau))$ for $t \ge 0$ \cite{csomos22aso}. By judiciously approximating $Q$ and the integrals involved in the Magnus
expansion, it is possible to devise a 2nd-order integrator with good preservation properties \cite{csomos20mti}. The analysis has been generalized in
\cite{csomos22aso} to more general classes of problems where $A$ depends in a more involved way of $\tau$-history $x_t(\theta) = x(t+\theta)$. 
Achieving higher orders of convergence with this procedure, however, is problematic and in fact it is left as open in \cite{csomos22aso}.	   

In the following, we show how the same procedure we have applied in the linear case can be generalized also in this setting, leading to higher order numerical
integrators. As before, we discretize the initial
function $\phi(t)$ with Chebyshev collocation points in the interval $[-\tau,0]$ and replace the abstract Cauchy problem with the 
$d(N+1)$-dimensional autonomous system
\begin{equation} \label{nl.1}
\begin{aligned}
 &  \dot{U}_N(t) = \mathcal{A}_N(U_N) \, U_N(t), \quad t \ge 0 \\
 & U_N(0) = \phi_N
\end{aligned} 
\end{equation}
where
\begin{equation} \label{nl.2}
  \mathcal{A}_N(U_N) = \left( \begin{array}{ccccc}
  			A((U_N)_{N+1}) & 0_d & \ldots & 0_d & 0 \\
			    &       & \frac{2}{\tau} [ \mathbb{D}^{(d+1,d(N+1))}] &  
			  \end{array} \right).  
\end{equation}
In other words, we replace in (\ref{lp.4}) $A(t)$ by the new $A$, and $B$ by the zero matrix.
Observe that, since the matrix $A$ in (\ref{eq.5.1}) only depends on $x(t-\tau)$, then the dependence of the matrix $\mathcal{A}_N$ in the
discretized system (\ref{nl.1}) comes only through the last $d$ components of the vector $U_N$, corresponding to the $(N+1)$-th collocation point. This is the
meaning of the notation $A((U_N)_{N+1})$ in (\ref{nl.2}).

Next, the nonlinear initial value problem (\ref{nl.2}) is solved with exponential integrators based on a generalized Magnus expansion. In fact, it is shown in \cite{casas06eme} how
methods up to order four in this class can be obtained for nonlinear equations of the form
\begin{equation} \label{nl.2b}
   \dot{y} = \hat{A}(t,y) \, y, \qquad y(0) = y_0.
\end{equation}   
As in the linear case, the starting point is to represent the solution in the form $y(t) = \exp(\Omega(t)) y_0$. Then, $\Omega(t)$ can be obtained by
Picard's iteration as
\begin{equation} \label{nl.3}
\begin{aligned}
  & \Omega^{[0]}(t) \equiv 0 \\
  & \Omega^{[m+1]}(t) = \int_0^t \sum_{k=0}^{\infty} \frac{B_k}{k!} \ad_{\Omega^{[m]}(s)}^k \hat{A} \big(s, \e^{\Omega^{[m]}(s)} y_0 \big) ds, \qquad m \ge 0
\end{aligned}
\end{equation}
so that $\lim_{m \rightarrow \infty} \Omega^{[m]}(t) = \Omega(t)$ for sufficiently small values of $t$. If terms up to $k=m-2$ are kept in (\ref{nl.3}), then
\[
 \Omega^{[m]}(t) = \Omega(t) + \mathcal{O}(t^{m+1}),
\]
so that by appropriate quadratures, it is then possible to get consistent approximations for the first values of $m$. In particular, one can construct the
following schemes of order 2 and 3, whereas the explicit algorithm for the method of order 4 can be found in \cite{casas06eme,hajiketabi20nib}.
For simplicity, we only consider the autonomous case in (\ref{nl.2b}).

\paragraph{Order 2.} Applying the trapezoidal rule to $\Omega^{[2]}$, one gets
\begin{equation} \label{nl.4}
\begin{aligned}
 & u = h \hat{A}(y_k) \\
 & v = \frac{1}{2} \left( u + h \hat{A}(\e^u y_k) \right) \\
 & y_{k+1} = \e^v y_k
\end{aligned}
\end{equation} 

\paragraph{Order 3.} Now one has to approximate the integrals appearing in $\Omega^{[1]}$, $\Omega^{[2]}$, and $\Omega^{[3]}$. The minimum number
of evaluations of $\hat{A}$ and commutators is achieved by the following sequence:
\begin{equation} \label{nl.5}
\begin{aligned}
 & Q_1 = h  \hat{A}(y_k) \\
 & Q_2 = h \hat{A}(\e^{\frac{1}{2} Q_1} y_k) - Q_1 \\
 & u_1 = \frac{1}{2} Q_1 + \frac{1}{4} Q_2 \\
 & u_2 = Q_1 + Q_2 \\
 & Q_3 = -u_2 + h \hat{A}(\e^{u_1} y_k) \\
 & Q_4 = -u_2 - Q_2 + h \hat{A}(\e^{u_2} y_k) \\
 & u_3 = u_2 + \frac{2}{3} Q_3 + \frac{1}{6} Q_4 - \frac{1}{6} [Q_1, Q_2] \\
 & y_{k+1} = \e^{u_3} u_k
\end{aligned}
\end{equation} 
Notice that these Magnus integrators require more computational effort than the corresponding schemes for the linear case. 

In the following we check the whole procedure on two additional numerical examples. As before, the different integrators are denoted by M2, M3 and M4,
respectively.

\begin{figure}[!ht] 
\begin{center}
\includegraphics[width=6.5cm]{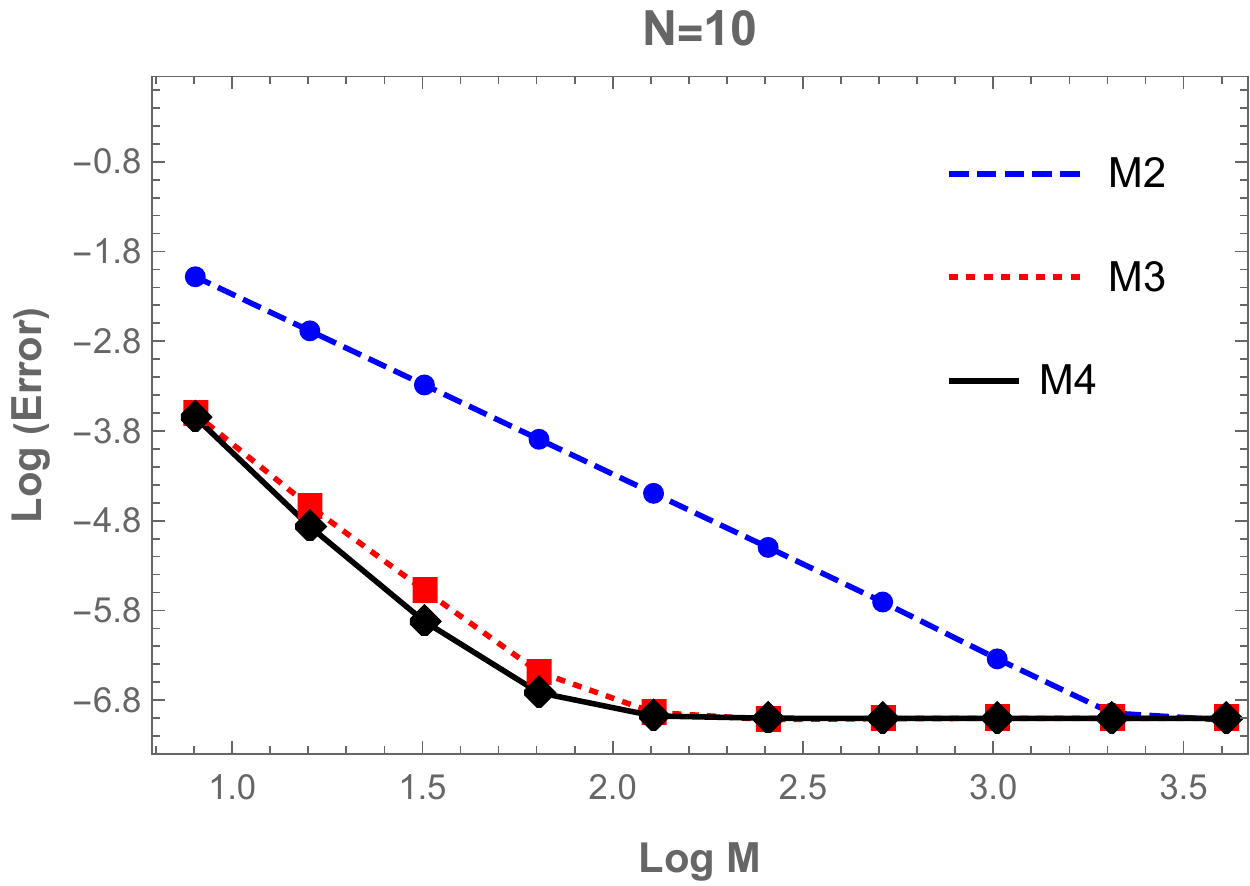} 
\includegraphics[width=6.5cm]{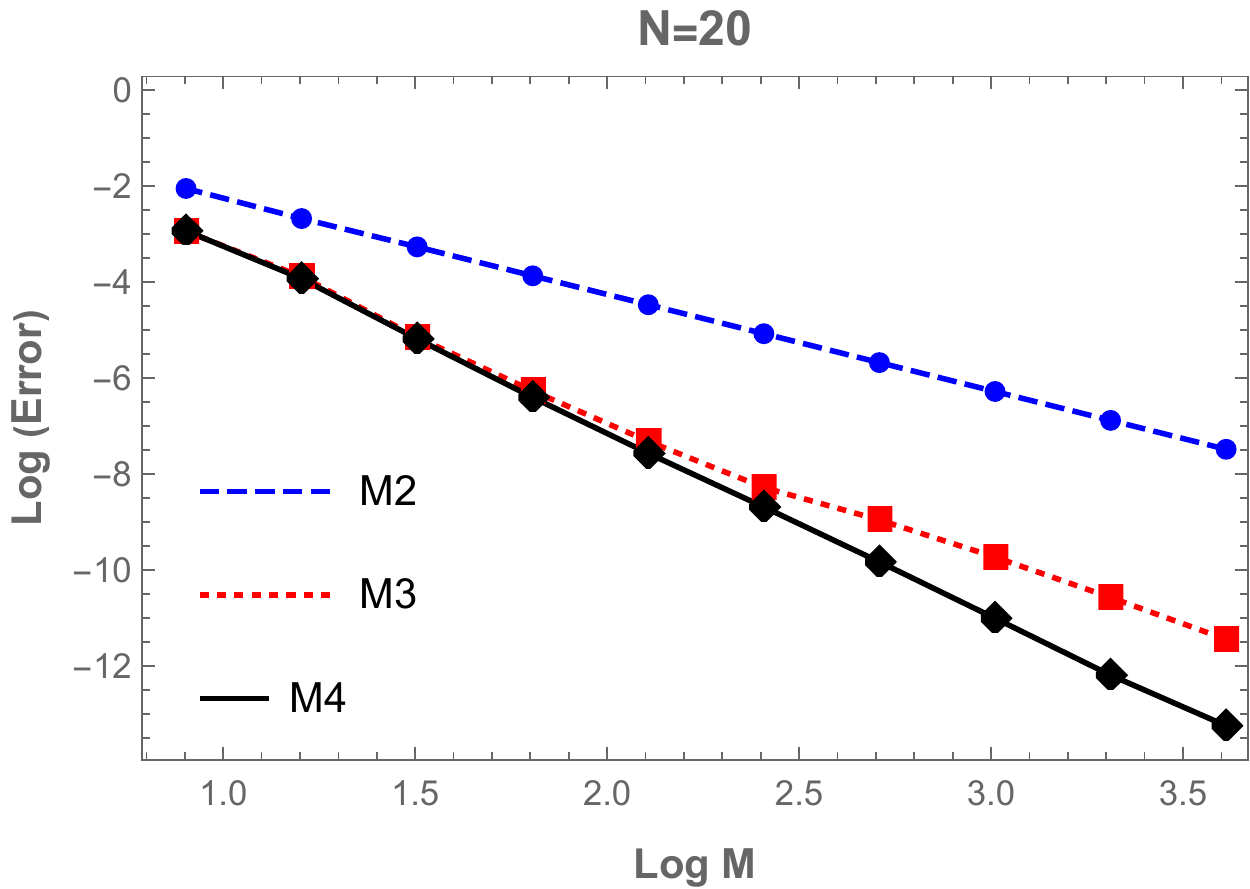} 
\end{center}
\caption{\label{figure5} \small Mean error in the solution obtained with Magnus integrators M2, M3 and M4 for the scalar nonlinear equation (\ref{1d.2}). Left: $N=10$. Right: $N=20$.}
\end{figure}

\paragraph{Example 3: a scalar nonlinear equation.} The scalar equation (\ref{1d.2}) constitutes a particularly simple example of a quasilinear
DDE of type (\ref{dde.2}) where the previous algorithm can be checked, since the exact solution is known explicitly, namely
$z(t) = \e^{\sin t}$. To do that, we integrate the resulting $N$-dimensional system (\ref{nl.1})-(\ref{nl.2}) from a discretization based on $N$ Chebyshev
collocation points (since $d=1$ here) until the final time $t_f = \tau = \frac{\pi}{2}$. As in Example 1, we compute the main error in the interval as a function
of $M$ for $N=10$ and $M=20$ to check that the Magnus integrators provide indeed the prescribed order. In this way we get Figure \ref{figure5}. Notice
that, whereas the order of M2 is clearly visible, the accuracy of M4 only manifests itself for a sufficiently large number of Chebyshev points and very small step
sizes. Otherwise, the results provided by M3 and M4 are quite similar.

\paragraph{Example 4: a SIR model with delay.}
In reference \cite{csomos20mti}, a 2nd-order numerical integrator also based on the Magnus expansion is proposed an tested on a widely
used delayed SIR model (see \cite{csomos20mti} and references therein), namely
\begin{equation} \label{sir.1}
\begin{aligned}
 & \dot{S}(t) = - \beta S(t) \, \frac{I(t-\tau)}{1+ \alpha I(t-\tau)} \\
 & \dot{I}(t) = \beta S(t) \,  \frac{I(t-\tau)}{1+ \alpha I(t-\tau)} - \gamma I(t) \\
 & \dot{R}(t) = \gamma I(t).
\end{aligned}
\end{equation}
It describes an epidemic model where the variation in the number of susceptible, $S(t)$, and infected, $I(t)$, individuals depend not only on the
actual values of $S(t)$ and $I(t)$ respectively, but also on how many infected individuals they interacted a latent period before, i.e., at $t-\tau$. In
(\ref{sir.1}), $\beta > 0$ and $\gamma > 0$ denote the infection and recovery rates, respectively, whereas $\alpha = 0$ if there is only
a bilinear incidence rate, and $\alpha = 1$ if the model incorporates a saturated incidence rate \cite{csomos20mti}. 

\begin{figure}[!ht] 
\begin{center}
\includegraphics[width=6.5cm]{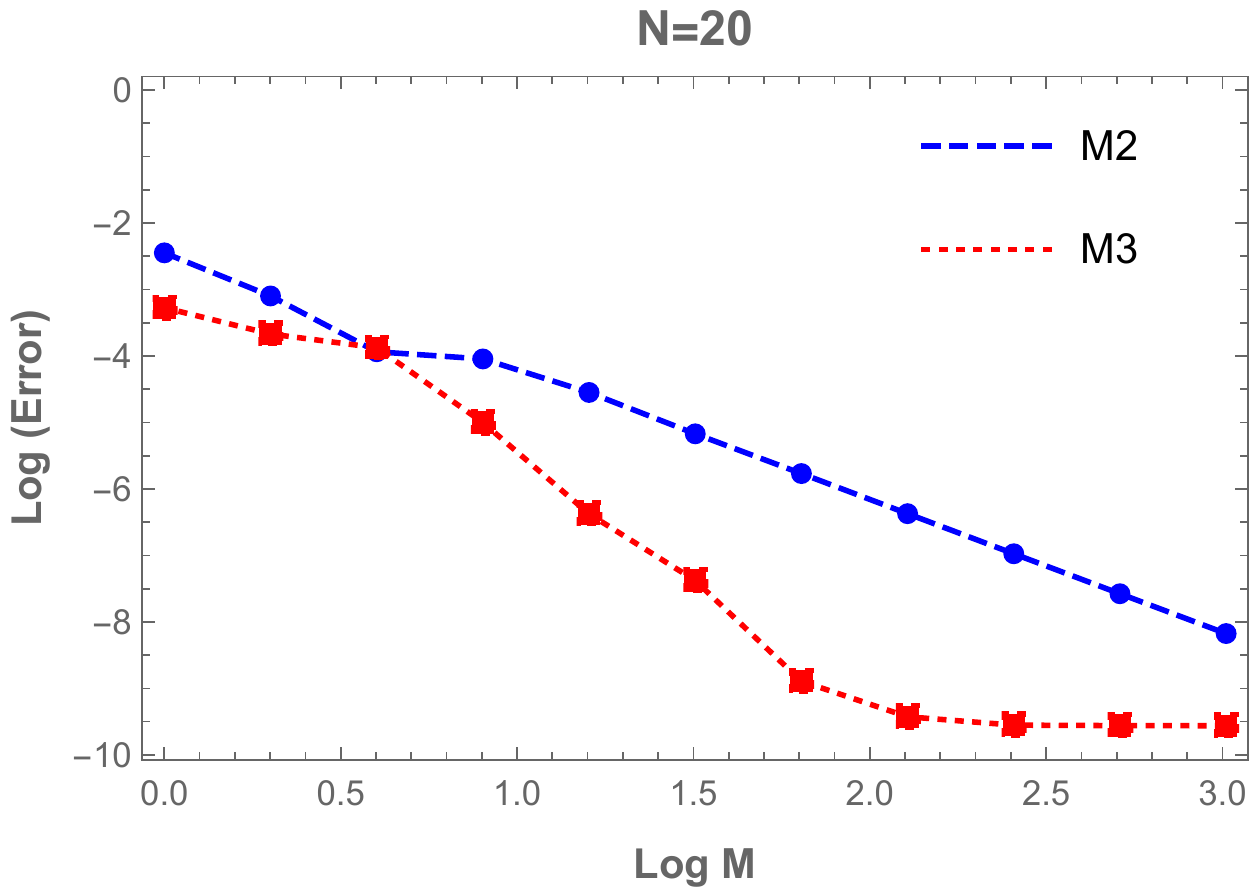} 
\includegraphics[width=6.5cm]{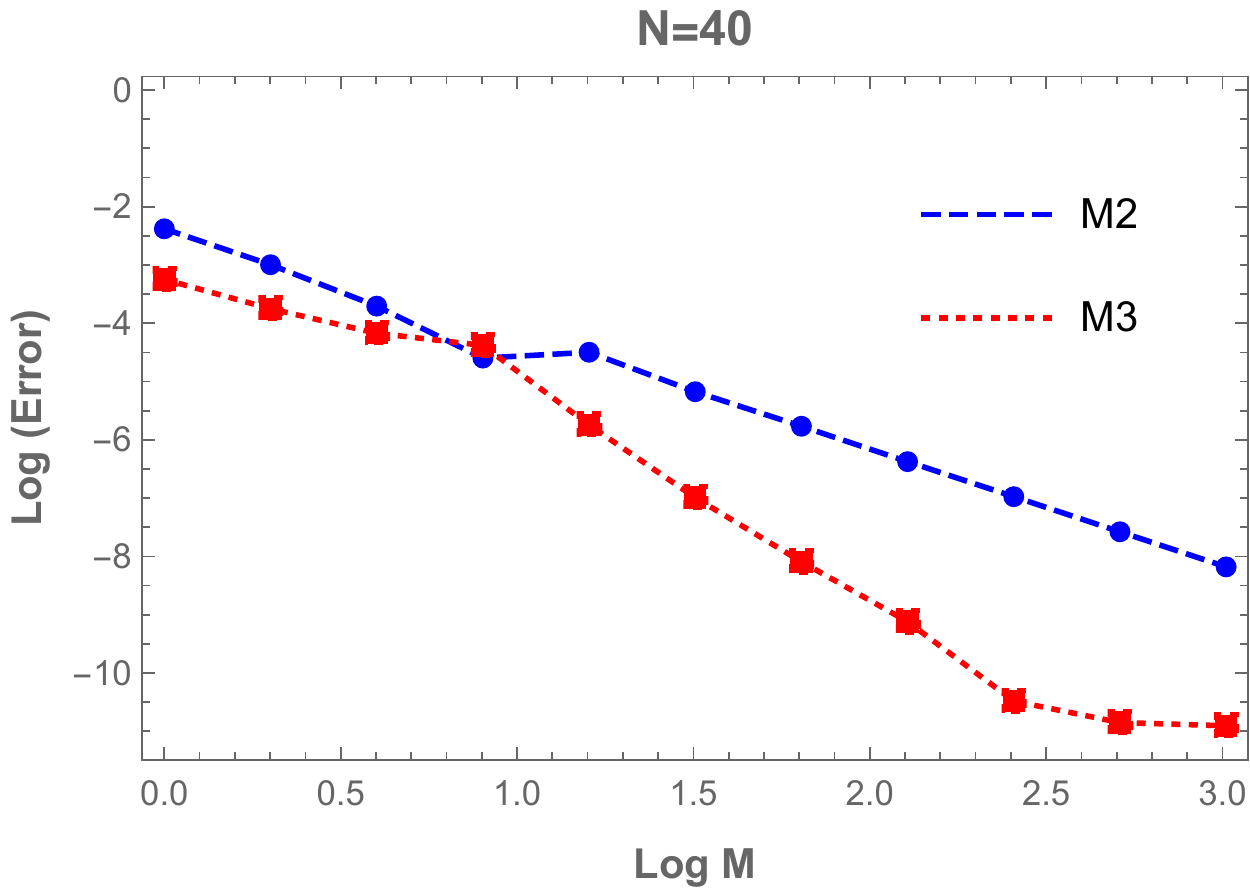} 
\end{center}
\caption{\label{figure6} \small Relative error in the solution obtained with M2 and M3 for the SIR model with delay (\ref{sir.1}). Left: $N=20$. Right: $N=40$.}
\end{figure}

Notice that system (\ref{sir.1}) can be expressed as eq. (\ref{dde.2}) with $x = (S, I, R)^T$ and
\begin{equation} \label{sir.2}
  A(x(t-\tau)) = \left( \begin{array}{crr}
  		-q(I(t-\tau)) & 0 & 0 \\
		q(I(t-\tau)) & -\gamma & 0 \\
		0 & \gamma & 0
	\end{array} \right),
\end{equation}
where $q(I(t- \tau)) = \frac{I(t-\tau)}{1+ \alpha I(t-\tau)}$. 
As in  \cite{csomos20mti}, we take
$\alpha = 0$, $\beta=1$, $\gamma = 1$, latent period $\tau = 1$, initial values $S_0 = 0.7$, $I_0 = 0.2$,
$R_0 = 0.1$ and initial function
\[
  \phi(t) = I_0 - \frac{1}{2} t,
\]  
for which the effect of the latent period is most evident. Then we integrate numerically with our algorithm until the final time $t_f = 4 \tau$ and compute
the relative error
\[
  \varepsilon_r = \frac{\|x(t_f) - x_n(t_f)\|}{\|x(t_f)\|},
\]  
where the `exact' solution $x(t_f)$ is obtained with the Matlab built-in function \texttt{dde23} with relative tolerance $10^{-12}$ and $x_n(t_f)$ refers
to the numerical solution computed with Magnus integrators (as appropriately extracted from the whole vector $U_N(t_f)$).
In view of the close similarity of the results achieved by M3 and M4 in the previous example, only M2 and M3 are tested here. The results we achieve are
shown in Figure \ref{figure6}.

Here again the order of each integrator is clearly visible, with M3 leading to high accuracy: with $h=0.01$ (corresponding to $M=100$) the relative error is
approximately $10^{-9}$. Now we get errors smaller than $10^{-2}$ with the largest possible value of the time step size, namely $h=\tau$.

\begin{figure}[!ht] 
\begin{center}
\includegraphics[scale=0.75]{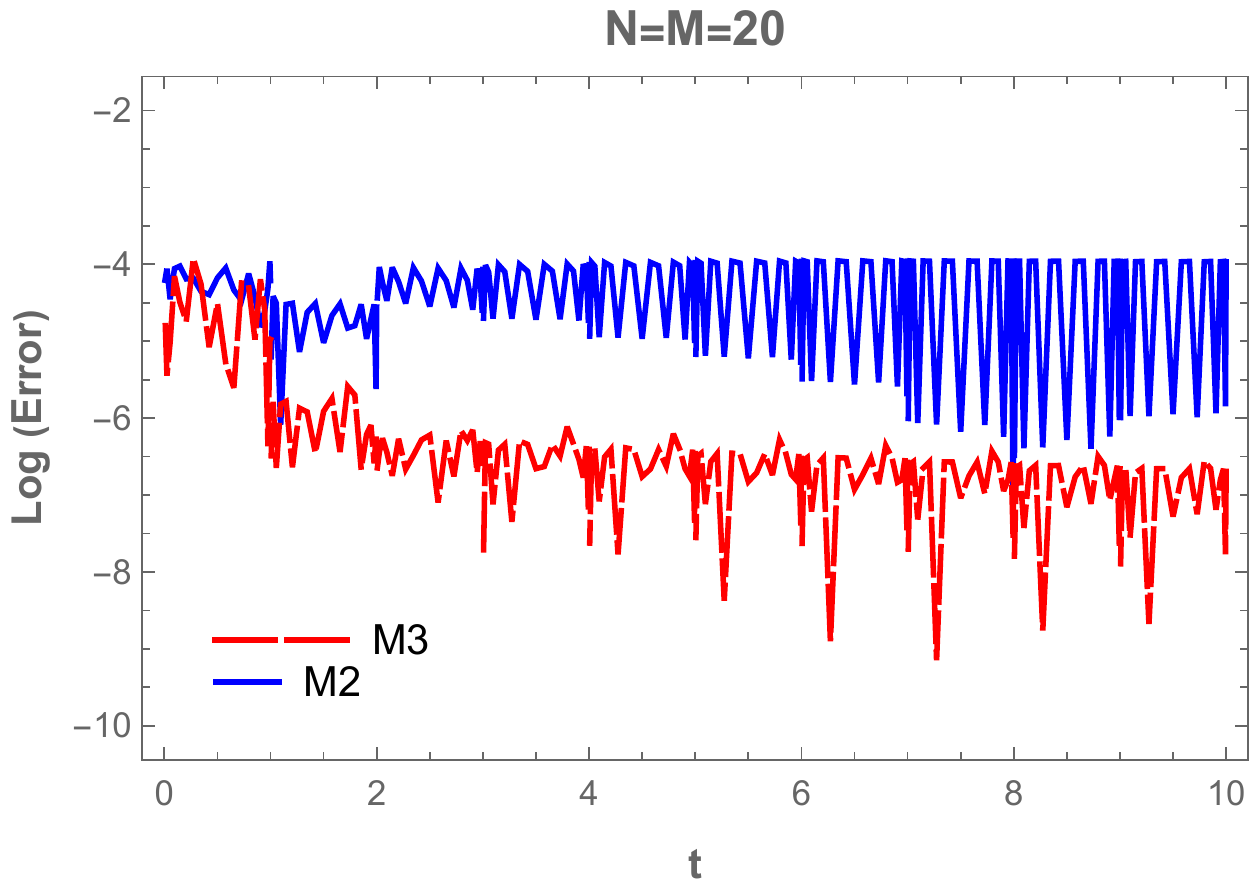} 
\end{center}
\caption{\label{figure7} \small Relative error as a function of time for the SIR model with delay (\ref{sir.1}) and $N=M=20$. 
Parameters $\alpha = 0$, $\beta=1$, $\gamma = 1$ and
initial function (\ref{infun2}).
Final time: $t_f = 10 \tau$.}
\end{figure}

Finally, in Figure \ref{figure7} we show the relative error as a function of time when the integration is carried out in the interval $[0, 10 \tau]$ with the same
values for the parameters,
initial function
\begin{equation} \label{infun2}
  \phi(t) = I_0 + \frac{1}{2} t,
\end{equation}  
$N=20$ collocation points and $M = h^{-1} = 20$. In this case, no secular growth in the error is observed, with M3 providing more accurate results.

\paragraph{Preservation of properties.} Since model (\ref{sir.1}) deals with a population, it is clear that $S(t)$, $I(t)$ and $R(t)$ are non-negative for 
all $t \ge 0$, and moreover $S(t) + I(t) + R(t)$ is constant. This feature is connected with the special structure the corresponding
matrix $A(x(t-\tau))$ possesses in this case: in fact, $A$ is a graph Laplacian, characterized by the following properties \cite{blanes22ppm}:
\begin{itemize}
 \item its elements $A_{k,l}$ verify that $A_{k,l} \ge 0$ for $k,l = 1, \ldots, d$, $k \ne l$, $A_{k,k} \le 0$ for $k=1,\ldots, d$;
 \item $\sum_{k=1}^d A_{k,l} = 0$ for $l=1,\ldots, d$.
\end{itemize} 
It is then natural to analyze whether the corresponding approximations 
obtained by our numerical algorithm also verify these properties.

\begin{figure}[!ht] 
\begin{center}
\includegraphics[scale=0.75]{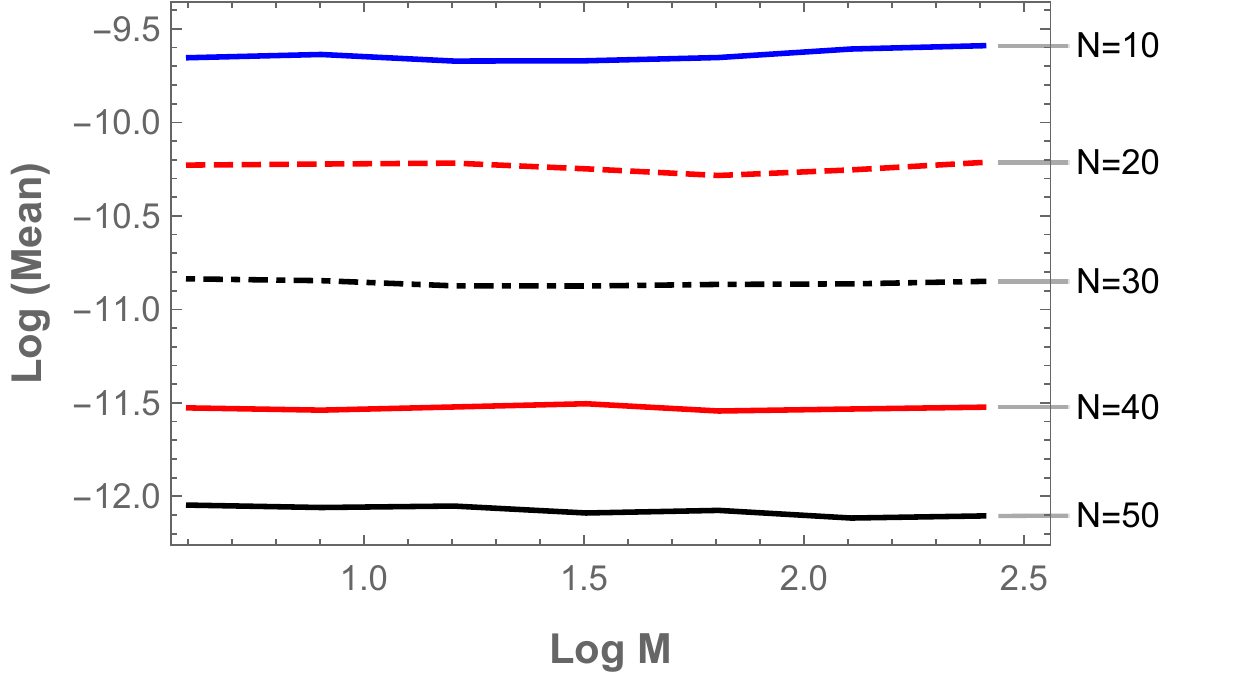} 
\end{center}
\caption{\label{figure8} \small Average value of the error in the total population (\ref{total1}) obtained with M3 for different values of $N$ and $M$ for the epidemic
model (\ref{sir.1}). The parameters and initial condition are the same as in Figure \ref{figure6}.}
\end{figure}

One can hardly expect that the matrix $\mathcal{A}_N(U_N)$ in (\ref{nl.1}) inherits the special structure $A$ may have  in general, due to the presence of 
the differentiation matrix $\mathbb{D}^{(d+1,d(N+1))}$. One should take into consideration, however, the following key observations:
\begin{itemize}
 \item the first $d$ rows of $\mathcal{A}_N(U_N)$ are zero, except for the first $d$ columns, where $A$ is placed; in other words,
 the first $d \times d$ block of $\mathcal{A}_N(U_N)$ is precisely $A$;
 \item the exponential matrix $\e^{\mathcal{A}_N}$ has the same structure as $\mathcal{A}_N$, and the first $d \times d$ block corresponds to $\e^{A}$;
 \item the computations involved in the Magnus integrators (\ref{nl.4}) and (\ref{nl.5}) do not alter this structure of the matrices; in other words, the matrices
 obtained at the intermediate stages M2 and M3 ($u, v, Q_1, Q_2, \ldots, u_3$) have the same basic structure as $\mathcal{A}_N$; 
 \item if all the elements of $w \in \mathbb{R}^d$ are non-negative, then all elements of $\e^u w$ are non-negative if $u$ is  graph Laplacian \cite{csomos20mti,blanes22ppm}.
\end{itemize}
In consequence, by construction, the previous Magnus integrators preserve the total population and positivity unconditionally at times $t=n \tau$, $n=1,2,\ldots$
\cite{blanes22ppm},
whereas at intermediate times the error is dictated by number of Chebyshev nodes used in the spectral discretization. This feature is clearly
illustrated in Figure
\ref{figure8}, where we show the average value of 
\begin{equation} \label{total1}
  |S(\tau + \theta_j) + I(\tau + \theta_j) + R(\tau + \theta_j) - 1|
\end{equation}
over the interval $[3 \tau, 4 \tau]$ obtained with the 3th-order Magnus integrator with different Chebyshev points and increasingly large numbers of subdivisions $M$.
Whereas at $t=\tau$, this difference is of the order of round-off, at intermediate steps it can be reduced significantly by increasing
the number $N$ of nodes. In all cases, the error does not grow with time.

\section{Concluding remarks}
\label{sec.6}

We have presented a numerical procedure to integrate non-autonomous linear delay differential equations consisting of two main components. 
First, the use
of an spectral discretization of the delayed part to transform the original problem into a  linear ordinary differential equation (ODE) whose dimension depends
on the number $N$ of collocation points: thus, if the original matrix system has dimension $d$, then the transformed ODE has dimension $d(N+1)$. 
Second, this system is then solved by applying numerical integrators based on the Magnus expansion. Thus, for a sufficiently large value of $N$, 
approximate solutions up to order 6 have been obtained for linear systems. In fact, if the step size is small enough, then the final accuracy is dictated by $N$. 
Even for large time steps, there are not stability issues related with the Magnus integrators.

The combination of both the spectral discretization and high order Magnus integrators 
leads to very accurate approximations to the characteristic multipliers and the solution of the problem. This can
be of particular interest when analyzing the stability of linear DDEs in practice.

The technique is thus different from the standard approach based on the use of a discrete method for ODEs (e.g., a Runge--Kutta scheme)
endowed with some interpolant, and it is
particularly flexible for this type of equations. In particular, it allows one to adjust the step size and even the order of the Magnus method along the
integration, according with a specified tolerance. Although high orders of convergence are observed in practice, establishing rigorously this property
is far from trivial, as the treatment done in \cite{csomos22aso} for a second order integrator shows.

Although only one delay has been considered here, the generalization to $k$ distinct delays is straightforward. We have also extended
the treatment to quasilinear delay equations of the form (\ref{dde.2}), this time in combination with integrators based on the Magnus expansion for nonlinear
equations. The treatment of the epidemic model with delay (\ref{sir.1}) shows that the algorithm preserves unconditionally important features of the
system, such as the total population and the positivity of the variables at times $t = n \tau$, whereas error at intermediate times is governed by the number of nodes $N$, which in any case remains constant.

The use of integrators based on the Magnus expansion involves the computation of the exponential of matrices of dimension $d(N+1) \times d(N+1)$,
and this is very often the most expensive part of the algorithm. In this respect, it is worth remarking that the number $N$ of collocation points is not
very large (in our examples we get excellent results already with $N=20$ or $N=30$), and that we use the technique presented in \cite{bader19ctm} to approximate
the exponential by conveniently chosen Taylor polynomials. In this way the computational cost is reduced with respect to standard algorithms based on
Pad\'e approximants.


\subsection*{Acknowledgements}
 This work has been supported by 
Ministerio de Ciencia e Innovaci\'on (Spain) through project PID2019-104927GB-C21, MCIN/AEI/10.13039/501100011033, and by Universitat Jaume I
through project UJI-B2019-17. AA is additionally funded by project PGC2018-094889-B-100, MCIN/AEI/10.13039/501100011033, ERDF (``A way of ma\-king Europe").


\end{document}